\documentclass[preprint,12pt,authoryear]{elsarticle}
\usepackage{amsmath,color}
\usepackage{amssymb}
\usepackage{amsfonts}
\usepackage{longtable}
\usepackage{tikz}
\usepackage{pgfplots}
\usepackage{graphicx}
\usepackage{csquotes}
\usepackage{dsfont}
\usepackage{pdflscape}
\usepackage{morefloats}
\usepackage{nccfoots}
\usepackage{enumitem} 
\usepackage[round]{natbib}
\usepackage{booktabs}
\usepackage{inputenc}
\usepackage{hyperref}

\newtheorem{theorem}{Theorem}[section]

\newtheorem{Lem}[theorem]{Lemma}
\newtheorem{Kor}[theorem]{Corollary}

\newtheorem{definition}[theorem]{Definition}
\newproof{pf}{Proof}
\newproof{pot}{Proof of Theorem \ref{thm:opd_dep_measure}}
\newproof{pot2}{Proof of Lemma \ref{lemma:AR(1)}}
\newproof{pot3}{Proof of Proposition \ref{prop:OPD_and_Kendall}}

\newtheorem{remark}{Remark}[section]
\newtheorem{example}{Example}[section]

\usepackage{amssymb}

\journal{Entropy}

\newcommand{\ignore}[1]{}

\makeatletter
\newcommand*{\defeq}{\mathrel{\rlap{%
                     \raisebox{0.3ex}{$\m@th\cdot$}}%
                     \raisebox{-0.3ex}{$\m@th\cdot$}}%
                     =}
\makeatother

\makeatletter
\newcommand*{\eqdef}{=\mathrel{\rlap{%
                     \raisebox{0.3ex}{$\m@th\cdot$}}%
                     \raisebox{-0.3ex}{$\m@th\cdot$}}%
                     }
\makeatother

\usepackage{listings}
\begin{document}

\begin{frontmatter}


 \title{Ordinal Pattern Dependence in the Context of Long-Range Dependence}
 \author{Ines N\"ußgen \corref{cor2}\fnref{label1}} \ead{nuessgen@mathematik.uni-siegen.de}
  \author{Alexander Schnurr \fnref{label1}} 
 \fntext[label1]{Siegen University, Emmy-Noether-Campus, Walter-Flex-Str. 3}
 \cortext[cor2]{Corresponding author}

\begin{abstract}
Ordinal pattern dependence is a multivariate dependence measure based on the co-movement of two time series. In strong connection to ordinal time series analysis, the ordinal information is taken into account to derive robust results on the dependence between the two processes. This article deals with ordinal pattern dependence for long-range dependent time series including mixed cases of short- and long-range dependence. We investigate the limit distributions for estimators of ordinal pattern dependence. In doing so we point out the differences that arise for the underlying time series having different dependence structures. Depending on these assumptions, central and non-central limit theorems are proven. The limit distributions for the latter ones can be included in the class of multivariate Rosenblatt processes. Finally, a simulation study is provided to illustrate our theoretical findings.
\end{abstract}


\begin{highlights}
\item Non-central limit theorems for estimators of ordinal pattern dependence are proven for pure long-range dependence 
\item An extension to limit theorems for mixed cases of short- and long-range dependent components in the multivariate processes is provided 
\item A simulation study illustrates limit distributions in the context of long-range dependence and the Rosenblatt distribution
\end{highlights}


\end{frontmatter}
\section{Introduction}

The origin of the concept of ordinal patterns is in the theory of dynamical systems. The idea is to consider the order of the values within a data vector instead of the full metrical information. The ordinal information is encoded as a permutation (cf. Section 3). Already in the first papers on the subject, the authors considered entropy concepts related to this ordinal structure (cf. \cite{bandt:pompe:2002}). There is an interesting relationship between these concepts and the well-known Komogorov-Sinai entropy (cf. \cite{keller:sinn:2010}, \cite{gutjahr:keller:2020}). Additionally, an ordinal version of the Feigenbaum diagram has been dealt with e.g. in  \cite{keller:sinn:2005}. 
In \cite{sinn:keller:2011}, ordinal patterns were used in order to estimate the Hurst parameter in long-range dependent time series. Hence, the concept made its way into the area of statistics. Instead of long patterns (or even letting the pattern length tend to infinity), rather short patterns have been considered in this new framework. 
Furthermore, ordinal patterns have been used in the context of ARMA processes (\cite{bandt:shiha:2007}) and change-point detection within one time series \cite{unakafov:2018}.
In \cite{schnurr:2014} ordinal patterns were used for the first time in order to analyze the dependence between two time series. Limit theorems for this new concept were proved in a short-range dependent framework in \cite{schnurr:dehling:2017}. Ordinal pattern dependence is a promising tool, which has already been used on financial, biological and hydrological data sets. Since in particular in this last context, the data sets are known to be long-range dependent. It is important to have limit theorems available also in this framework. We close this gap in the present article. 

All of the results presented in this article have been established in the PhD-thesis of I. Nüßgen written under the supervision of A. Schnurr. 

The article is structured as follows: in the subsequent section we provide the reader with the mathematical framework. The focus is on (multivariate) long-range dependence. In Section 3 we recall the concept of ordinal pattern dependence and prove our main results. We present a simulation study in Section 4 and close the paper by a short outlook in Section 5. 

\section{Mathematical framework}
\noindent We consider a stationary $d$-dimensional Gaussian time series $\left(Y_j\right)_{j\in\mathbb{Z}}$ (for $d\in\mathbb{N}$), with
\begin{align}
 Y_j:=\left(Y_j^{(1)},\ldots,Y_j^{(d)}\right)^t\label{eq:process}
\end{align}
such that $\mathbb{E}\left(Y_j^{(p)}\right)=0$ and $\mathbb{E}\left(\left(Y_j^{(p)}\right)^2\right)=1$ for all $j\in \mathbb{Z}$ and $p=1,\ldots,d$. Furthermore, we require the cross-correlation function to fulfill $\left|r^{(p,q)}(k)\right|<1$ for $p,q=1,\ldots,d$ and $k\geq 1$, where the componentwise cross-correlation functions $r^{(p,q)}(k)$ are given by $r^{(p,q)}(k)=\mathbb{E}\left(Y_j^{(p)}Y_{j+k}^{(q)}\right)$ for each $p,q=1,\ldots,d$ and $k\in\mathbb{Z}$. For each random vector $Y_j$ we denote the covariance matrix by $\Sigma_{d}$, since it is independent of $j$ due to stationarity. Therefore, we have $\Sigma_d=\left(r^{(p,q)}(0)\right)_{p,q=1,\ldots, d}$.  \newline\newline
We specify the dependence structure of $\left(Y_j\right)_{j\in\mathbb{Z}}$ and turn to long-range dependence: we assume that for the cross-correlation function $r^{(p,q)}(k)$ for each $p,q=1,\ldots, d$, it holds that
\begin{align}
r^{(p,q)}(k)= L_{p,q}(k) k^{d_p+d_q-1}, \label{multivariateLRDconditionII}
\end{align}
with $ L_{p,q}(k) \rightarrow L_{p,q}$ $(k\rightarrow\infty)$ for finite constants $L_{p,q}\in[0,\infty)$ with $L_{p,p}\neq 0$, where the matrix $L=\left(L_{p,q}\right)_{p,q=1,\ldots,d}$ has full rank, is symmetric and positive definite. Furthermore, the parameters $d_p,d_q\in \left(0,\frac{1}{2}\right)$ are called long-range dependence parameters.
Therefore, $\left(Y_j\right)_{j\in\mathbb{Z}}$ is multivariate long-range dependent in the sense of \cite{kechagias:pipiras:2015}, Definition 2.1.\newline\newline
The processes we want to consider have a particular structure, namely for $h\in\mathbb{N}$, we obtain for fixed $j\in\mathbb{Z}$:
\begin{align}
Y_{j,h}:=\left(Y_j^{(1)},\ldots,Y_{j+h-1}^{(1)},Y_j^{(2)},\ldots,Y_{j+h-1}^{(2)},\ldots,Y_j^{(d)},\ldots,Y_{j+h-1}^{(d)}\right)^t\in\mathbb{R}^{dh}.\label{eq:multivGP}
\end{align}
The following relation holds between the \emph{extendend process} $\left(Y_{j,h}\right)_{j\in\mathbb{Z}}$ and the primarily regarded process $\left(Y_j\right)_{j\in\mathbb{Z}}$. For all $k=1,\ldots, dh$, $j\in\mathbb{Z}$ we have
\begin{align}
Y_{j,h}^{(k)}=Y_{j+(k\mod h)-1}^{\left\lfloor \frac{k-1}{h}\right\rfloor +1}, \label{eq: relation extended process to original}
\end{align}
where $\lfloor x \rfloor=\max\{k\in\mathbb{Z}:k\leq x\}$.
Note that the process $\left(Y_{j,h}\right)_{j\in\mathbb{Z}}$ is still a centered Gaussian process since all finite-dimensional marginals of $\left(Y_j\right)_{j\in\mathbb{Z}}$ follow a normal distribution. Stationarity is preserved, too, since for all $p,q=1,\ldots,dh$, $p\leq q$ and $k\in\mathbb{Z}$ the cross-correlation function $r^{(p,q,h)}(k)$ of the process $\left(Y_{j,h}\right)_{j\in\mathbb{Z}}$ is given by \newline
\begin{align}
r^{(p,q,h)}(k)&=\mathbb{E}\left(Y_{j,h}^{(p)}Y_{j+k,h}^{(q)}\right)\notag\\
&=\mathbb{E}\left(Y_{j+(p\mod h)-1}^{\left\lfloor \frac{p-1}{h}\right\rfloor +1}Y_{j+k+(q\mod h)-1}^{\left\lfloor \frac{q-1}{h}\right\rfloor +1}\right)\notag\\
&=r^{(\left\lfloor \frac{p-1}{h}\right\rfloor +1,\left\lfloor \frac{q-1}{h}\right\rfloor +1)}(k+((q-p)\mod h))\label{eq:crosscorrelationextendedprocess}
\end{align} 
and the last line does not depend on $j$. The covariance matrix $\Sigma_{d,h}$ of $Y_{j,h}$ has the following structure:
\begin{align*}
\left(\Sigma_{d,h}\right)_{p,q=1,\ldots,d,\atop p\leq q}&=\left(r^{(p,q,h)}(0)\right)_{p,q=1,\ldots,dh,\atop p\leq q,}, \\
\left(\Sigma_{d,h}\right)_{p,q=1,\ldots,d,\atop p>q}&=\left(r^{(q,p,h)}(0)\right)_{p,q=1,\ldots,dh,\atop q<p}.
\end{align*}
Hence, we arrive at
\begin{align}
\Sigma_{d,h}=\left(\Sigma_{h}^{(p,q)}\right)_{1\leq p,q \leq d}, \label{eq: covariance matrix extended process}
\end{align}
where $\Sigma_h^{(p,q)}=\mathbb{E}\left(\left(Y_1^{(p)},\ldots,Y_h^{(p)} \right)^t\left(Y_1^{(q)},\ldots,Y_h^{(q)} \right)\right)=\left(r^{(p,q)}(i-k)\right)_{1\leq i,k\leq h}$, $p,q=1,\ldots,d$. Note that $\Sigma_h^{(p,q)}\in\mathbb{R}^{h\times h}$ and $r^{(p,q)}(k)=r^{(q,p)}(-k)$, $k\in\mathbb{Z}$ since we are studying cross-correlation functions.\newline
So finally we have to show that based on the assumptions on $\left(Y_j\right)_{j\in\mathbb{Z}}$ the extended process is  still long-range dependent. \newline
Hence, we have to consider the cross-correlations again
\begin{align}
r^{(p,q,h)}(k)&=r^{(\left\lfloor \frac{p-1}{h}\right\rfloor +1,\left\lfloor \frac{q-1}{h}\right\rfloor +1)}(k+((q-p)\mod h))\notag\\
&=r^{(p^*,q^*)}(k+m^*)\notag\\
&\simeq r^{(p^*,q^*)}(k)~~(k\rightarrow\infty),\label{eq:crosscorrelationsjh}
\end{align}
since $p^*,q^*\in\{1,\ldots,d\}$ and $m^*\in\{0,\ldots,h-1\}$, with $p^*:=\left\lfloor \frac{p-1}{h}\right\rfloor +1$, $q^*:=\left\lfloor \frac{q-1}{h}\right\rfloor +1$ and $m^*=(q-p)\mod h$. \newline
Let us remark that $a_k\simeq b_k \Leftrightarrow \lim_{k\rightarrow\infty} \frac{a_k}{b_k}=1$.\newline
\newline
Therefore, we are still dealing with a multivariate long-range dependent Gaussian process. We see in the proofs of the following limit theorems that the crucial parameters that determine the asymptotic distribution are the long-range dependence parameters $d_p$, $p=1,\ldots,d$ of the original process $\left(Y_j\right)_{j\in\mathbb{Z}}$ and therefore, we omit the detailed description of the parameters $d_{p^*}$ here. \newline
It is important to remark that the extended process $\left(Y_{j,h}\right)_{j\in\mathbb{Z}}$ is also long-range dependent in the sense of \cite{arcones:1994}, p. 2259, since
\begin{align}
\lim_{k\rightarrow\infty} \frac{k^{D} r^{(p,q,h)}(k)}{L(k)}&=\lim_{k\rightarrow\infty} \frac{k^{D} r^{(p^*,q^*)}(k) }{L(k)}\notag\\
&=\lim_{k\rightarrow\infty} \frac{k^{D} L_{p^*,q^*}k^{d_{p^*}+d_{q^*}-1} }{L(k)}\notag\\
&=:b_{p^*,q^*},\label{vergleichLRDconditions}
\end{align}
with 
\begin{align}
D:=\min\limits_{p^*\in\{1,\ldots,d\}}\{1-2d_{p^*}\}\in(0,1) \label{eq: arcones alpha}
\end{align} and $L(k)$ can be chosen as any constant $L_{p,q}$ that is not equal to zero, so for simplicity we assume without loss of generality $L_{1,1}\neq 0$ and, therefore, $L(k)=L_{1,1}$, since the condition in \cite{arcones:1994} only requires convergence to a finite constant $b_{p^*,q^*}$. Hence, we may apply the results in \cite{arcones:1994} in the subsequent results.\newline
We define the following set, which is needed in the proofs of the theorems of this section.
\newline
\begin{align}
P^*:=\{p\in\{1,\ldots,d\}:~d_p\geq d_q, \text{ for all } q\in\{1,\ldots,d\}\} \label{eq: set1}
\end{align}
and denote the corresponding long-range dependence parameter to each $p\in P^*$ by
\begin{align*}
d^*:=d_p,\quad p\in P^*.
\end{align*}\newline
We briefly recall the concept of Hermite polynomials as they play a crucial role in determining the limit distribution of functionals of multivariate Gaussian processes.
\begin{definition}(Hermite polynomial, \cite{beran:feng:ghosh:kulik:2013}, Definition 3.1)\newline
The $j$-th Hermite polynomial $H_j(x)$, $j=0,1,\ldots$, is defined as
\begin{align*}
H_j(x):=(-1)^j\exp\left(\frac{x^2}{2}\right)\frac{\mathrm{d}^j}{\mathrm{d}x^j}\exp\left(-\frac{x^2}{2}\right).
\end{align*}
\end{definition}
Their multivariate extension is given by the subsequent definition.
\begin{definition}(Multivariate Hermite polynomial, \cite{beran:feng:ghosh:kulik:2013}, p. 122)\newline
Let $d\in\mathbb{N}$. We define as $d$-dimensional Hermite polynomial
\begin{align*}
H_k(x):=H_{k_1,\ldots,k_d}(x):=H_{k_1,\ldots,k_d}\left(x_1,\ldots,x_d\right)=\prod_{j=1}^d H_{k_j}\left(x_j\right),
\end{align*}
with $k=\left(k_1,\ldots,k_d\right)\in\mathbb{N}_0^d\setminus\{(0,\ldots,0)\}$.
\end{definition}
Let us remark that the case $k=(0,\ldots,0)$ is excluded here due to the assumption $\mathbb{E}\left(f(X)\right)=0$.\newline
Analogously to the univariate case, the family of multivariate Hermite polynomials \newline
$\left\{H_{k_1,\ldots,k_d},~k_1,\ldots,k_d\in\mathbb{N}\right\}$ forms an orthogonal basis of $L^2\left(\mathbb{R}^d,\varphi_{I_d}\right)$, which is defined as 
\begin{align*}
L^2\left(\mathbb{R}^d,\varphi_{I_d}\right):=\left\{f:\mathbb{R}^d\rightarrow\mathbb{R},~\int_{\mathbb{R}^d} f^2\left(x_1,\ldots,x_d\right)\varphi\left(x_1\right)\ldots\varphi\left(x_d\right)\mathrm{d}x_d\ldots\mathrm{d}x_1<\infty\right\}.
\end{align*}
The parameter $\varphi_{I_d}$ denotes the density of the $d$-dimensional standard normal distribution, which is already divided into the product of the univariate densities $\varphi$ in the formula above.\newline
We denote the Hermite coefficients by
\begin{align*}
C(f,X,k):=C\left(f,I_d,k\right):=\langle f,H_k\rangle=\mathbb{E}\left(f(X)H_k(X)\right).
\end{align*}
The Hermite rank $m\left(f,I_d\right)$ of $f$ with respect to the distribution $\mathcal{N}\left(0,I_d\right)$ is defined as the largest integer $m$, such that
\begin{align*}
\mathbb{E}\left(f(X)\prod_{j=1}^d H_{k_j}\left(X^{(j)}\right)\right)=0 \text{   for all } 0<k_1+\ldots k_d<m.
\end{align*}
Having these preparatory results in mind, we derive a the multivariate Hermite expansion given by
\begin{align}
f(X)-\mathbb{E}f(X)=\sum_{k_1+\ldots+k_d\geq m\left(f,I_d\right)} \frac{C(f,X,k)}{k_1!\ldots k_d!} \prod_{j=1}^d H_{k_j}\left(X^{(j)}\right).  \label{eq:multivariateHermiteexpansion}
\end{align}
We focus on limit theorems for functionals with Hermite rank $2$. First, we introduce the matrix-valued Rosenblatt process. This plays a cruicial role in the asymptotics of functionals with Hermite rank $2$ applied to multivariate long-range dependent Gaussian processes. 
We begin with the definition of a multivariate Hermitian-Gaussian random measure $\tilde{B}(\mathrm{d}\lambda)$ with independent entries given by
\begin{align}
\tilde{B}(\mathrm{d}\lambda)=\left(\tilde{B}^{(1)}(\mathrm{d}\lambda),\ldots,\tilde{B}^{(d)}(\mathrm{d}\lambda)\right)^t,\label{eq:HGrandommeasuremultivariate}
\end{align}
where $\tilde{B}^{(p)}(\mathrm{d}\lambda)$ is a univariate Hermitian-Gaussian random measure as defined in  \cite{pipiras:taqqu:2017}, Definition B.1.3. The multivariate Hermitian-Gaussian random measure $\tilde{B}(\mathrm{d}\lambda)$ satisfies
\begin{align*}
\mathbb{E}\left(\tilde{B}(\mathrm{d}\lambda)\right)&=0,\\
\mathbb{E}\left(\tilde{B}(\mathrm{d}\lambda)\tilde{B}(\mathrm{d}\lambda)^*\right)&=I_d~\mathrm{d}\lambda\\
\intertext{and}
\mathbb{E}\left(\tilde{B}^{(p)}(\mathrm{d}\lambda_1)\overline{\tilde{B}^{(q)}(\mathrm{d}\lambda_2)}\right)&=0,\quad \left|\lambda_1\right|\neq\left|\lambda_2\right|,\quad p,q=1,\ldots,d,
\end{align*}
where $\tilde{B}(\mathrm{d}\lambda)^*=\left(\overline{B^{(1)}\left(\mathrm{d}\lambda\right)},\ldots,\overline{B^{(d)}(\mathrm{d}\lambda)}\right)$ denotes the Hermitian transpose of $\tilde{B}(\mathrm{d}\lambda)$. Thus, following \cite{arcones:1994}, Theorem 6, we can state the spectral representation of the matrix-valued Rosenblatt process $Z_{2,H}(t)$, $t\in[0,1]$ as 
\begin{align*}
Z_{2,H}(t)=\left(Z_{2,H}^{(p,q)}(t)\right)_{p,q=1,\ldots,d} 
\end{align*}
where each entry of the matrix is given by
\begin{align*}
Z_{2,H}^{(p,q)}(t)=\int_{\mathbb{R}^2}^{\prime\prime} \frac{\exp\left(it\left(\lambda_1+\lambda_2\right)\right)-1}{i\left(\lambda_1+\lambda_2\right)}\tilde{B}^{(p)}\left(\mathrm{d}\lambda_1\right)\tilde{B}^{(q)}\left(\mathrm{d}\lambda_2\right).
\end{align*}~\quad
\noindent The double prime in $\int_{\mathbb{R}^2}^{\prime\prime}$ excludes the diagonals $\left|\lambda_i\right|=\left|\lambda_j\right|$, $i\neq j$ in the integration. For details on multiple Wiener-It\^o integrals see \cite{major:2019}.\newline\newline
The following results are taken from \cite{nuessgen:2021}, Sec. 3.2. The corresponding proofs are outsourced to the Technical Appendix.
\begin{theorem}\label{th: Hermite 2}
Let $\left(Y_j\right)_{j\in\mathbb{Z}}$ be a stationary Gaussian process as defined in \eqref{eq:process} that fulfills \eqref{multivariateLRDconditionII} for $d_p\in \left(\frac{1}{4},\frac{1}{2}\right)$, $p=1\ldots,d$. For $h\in\mathbb{N}$ we fix 
\begin{align*}
Y_{j,h}:=\left( Y_j^{(1)},\ldots,Y_{j+h-1}^{(1)},\ldots,Y_j^{(d)},\ldots,Y_{j+h-1}^{(d)}  \right)^t\in\mathbb{R}^{dh}
\end{align*}
with $Y_{j,h}\sim\mathcal{N}\left(0,\Sigma_{d,h}\right)$ and $\Sigma_{d,h}$ as described in \eqref{eq: covariance matrix extended process}. Let $f:\mathbb{R}^{dh}\rightarrow\mathbb{R}$ be a function with Hermite rank $2$ such that the set of discontinuity points $D_f$ is a Null set with respect to the $dh$-dimensional Lebesgue measure. Furthermore, we assume $f$ fulfills $\mathbb{E}\left( f^2\left( Y_{j,h}  \right)\right)<\infty$. Then,
\begin{align}
n^{-2d^*}(C_2)^{-\frac{1}{2}}&\sum_{j=1}^n  \left(f\left(Y_j^{(1)},\ldots,Y_{j+h-1}^{(d)}\right)-\mathbb{E}\left(f\left(Y_j^{(1)},\ldots,Y_{j+h-1}^{(d)}\right)\right)\right)\notag\\&\xrightarrow{\mathcal{D}} \sum_{p,q\in P^*}\tilde{\alpha}^{(p,q)} Z^{(p,q)}_{2,d^*+1/2}(1),\label{aussageTheoremHR2}
\end{align}
where
\begin{align*}
Z^{(p,q)}_{2,d^*+1/2}(1)=K_{p,q}\left(d^*\right)\int_{\mathbb{R}^2}^{\prime\prime} \frac{\exp\left(i\left(\lambda_1+\lambda_2\right)\right)-1}{i\left(\lambda_1+\lambda_2\right)}\left|\lambda_1\lambda_2  \right|^{-d^*} \tilde{B}_L^{(p)}\left(\mathrm{d}\lambda_1\right)\tilde{B}_L^{(q)}\left( \mathrm{d}\lambda_2\right).
\end{align*}
The matrix $K\left(d^*\right)$ is a normalizing constant, see \cite{nuessgen:2021}, Corollary 3.6. Moreover, $\tilde{B}_L(\mathrm{d}\lambda)$ is a multivariate Hermitian-Gaussian random measure with $\mathbb{E}\left(B_L(\mathrm{d}\lambda)B_L(\mathrm{d}\lambda)^*\right)=L~\mathrm{d}\lambda$ and $L$ as defined in \eqref{multivariateLRDconditionII}.
Furthermore, $C_2:=\frac{1}{2d^*\left(4d^*-1\right)}$ is a normalizing constant and
\begin{align*}
\tilde{\alpha}^{(p,q)}:=\sum\limits_{i,k=1}^h \alpha_{i,k}^{(p,q)}
\end{align*}
 where $\alpha^{(p,q)}_{i,k}=\alpha_{i+(p-1)h,k+(q-1)h}$ for each $p,q\in P^*$ and $i,k=1,\ldots,h$ and 
  \begin{align*}
  \left(\alpha_{i,k}\right)_{1\leq i,k\leq dh}=\Sigma_{d,h}^{-1}C\Sigma_{d,h}^{-1}
  \end{align*}
where $C$ denotes the matrix of second order Hermite coefficients, given by
\begin{align*}
C=\left(c_{i,k}\right)_{1\leq i,k\leq dh}=\mathbb{E}\left(Y_{1,h}\left(f\left(Y_{1,h}\right)-\mathbb{E} \left(f\left(Y_{1,h}\right)\right)\right) Y_{1,h}^t\right).
\end{align*}
\end{theorem}
It is possible to soften the assumptions in Theorem \ref{th: Hermite 2} to allow for mixed cases of short- and long-range dependence.
\begin{Kor}\label{cor:LTHR2} Instead of demanding in the assumptions of Theorem \ref{th: Hermite 2} that \eqref{multivariateLRDconditionII} holds for $\left(Y_{j}\right)_{j\in\mathbb{Z}}$ with the addition that for all $p=1,\ldots,d$ we have $d_p\in\left(\frac{1}{4},\frac{1}{2}\right)$, we may use the following condition:\newline
We assume that
\begin{align*}
r^{(p,q)}(k)=k^{d_p+d_q-1}L_{p,q}(k)\quad(k\rightarrow\infty)
\end{align*}
with $L_{p,q}(k)$ as given in \eqref{multivariateLRDconditionII}, but we do no longer assume $d_p\in\left(\frac{1}{4},\frac{1}{2}\right)$ for all $p=1,\ldots,d$ but soften the assumption to $d^*\in\left(\frac{1}{4},\frac{1}{2}\right)$ and for $d_p\neq d^*$, $p=1,\ldots,d$ we allow for $d_p\in\left(-\infty,0\right)\cup\left(0,\frac{1}{4}\right]$. 
Then, the statement of Theorem \ref{th: Hermite 2} remains valid.
\end{Kor}
However, with a mild technical assumption on the covariances of the one-dimensional marginal Gaussian processes that is often fulfilled in applications, there is another way of normalizing the partial sum on the right-hand side in Theorem \ref{th: Hermite 2}, this time explicity for the case $\#P^*=2$ and $h\in\mathbb{N}$, such that the limit can be expressed in terms of two standard Rosenblatt random variables. This yields the possibility to further study the dependence structure between these two random variables. In the following theorem we assume $\#P^*=d=2$ for the reader's convenience.
\begin{theorem}\label{th: HR2alternative}
Under the same assumptions as in Theorem \ref{th: Hermite 2} with $\#P^*=d=2$ and $d^*\in \left(\frac{1}{4},\frac{1}{2}\right)$ and the additional condition that $r^{(1,1)}(l)=r^{(2,2)}(l)$, for $l=0,\ldots,h-1$, and $L_{1,1}+L_{2,2}\neq L_{1,2}+L_{2,1}$, it holds that
\begin{align*}
n^{-2d^*}(C_2)^{-\frac{1}{2}}&\sum_{j=1}^n  \left(f\left(Y_j^{(1)},\ldots,Y_{j+h-1}^{(d)}\right)-\mathbb{E}f\left(Y_j^{(1)},\ldots,Y_{j+h-1}^{(d)}\right)\right)\\
&\xrightarrow{\mathcal{D}} \left(\tilde{\alpha}^{(1,1)}-\tilde{\alpha}^{(1,2)}\right) \frac{L_{2,2}-L_{2,1}-L_{1,2}+L_{1,1}}{2}Z^*_{2,d^*+1/2}(1)\\
&\quad\quad+\left(\tilde{\alpha}^{(1,1)}+\tilde{\alpha}^{(1,2)}\right) \frac{L_{2,2}+L_{2,1}+L_{1,2}+L_{1,1}}{2}Z^{**}_{2,d^*+1/2}(1)
\end{align*}
with $C_2:=\frac{1}{2d^*\left(4d^*-1\right)}$ being the same normalizing factor as in Theorem \ref{th: Hermite 2}, $\left(\alpha_{i,k}\right)_{1\leq i,k\leq dh}=\Sigma_{d,h}^{-1}C\Sigma_{d,h}^{-1}$ and $C=\left(c_{i,k}\right)_{1\leq i,k\leq dh}=\mathbb{E}\left(Y_{1,h}\left(f\left(Y_{1,h}\right)-\mathbb{E} f\left(Y_{1,h}\right)\right) Y_{1,h}^t\right)$.
Note that $Z^{*}_{2,d^*+1/2}(1)$ and $Z^{**}_{2,d^*+1/2}(1)$ are both standard Rosenblatt random variables whose covariance is given by
\begin{align}
Cov\left(Z^{*}_{2,d^*+1/2}(1),Z^{**}_{2,d^*+1/2}(1)\right)=\frac{\left(L_{2,2}-L_{1,1}\right)^2  }{\left(L_{1,1}+L_{2,2}\right)^2-\left(L_{1,2}+L_{2,1}\right)^2}.\label{eq:covarianceHR2alternative}
\end{align} 
\end{theorem}
\section{Ordinal pattern dependence}\label{sec: OPD}
Ordinal pattern dependence is a multivariate dependence measure that compares the co-movement of two time series based on the ordinal information. First introduced in \cite{schnurr:2014} to analyze financial time series, a mathematical framework including structural breaks and limit theorems for functionals of absolutely regular processes has been built in \cite{schnurr:dehling:2017}.\newline\newline
We start with the definition of an ordinal pattern and the basic mathematical framework that we need to build up the ordinal model.\newline
Let $S_h$ denote the set of permutations in $\{0,\ldots,h\}$, $h\in\mathbb{N}_0$ that we express as $(h+1)$-dimensional tuples, assuring that each tuple contains each of the numbers above exactly once. In mathematical terms, this yields
\begin{align*}
S_h=\left\{\pi\in \mathbb{N}_0^{h+1}:\quad 0\leq\pi_i\leq h, \text{ and } \pi_i\neq\pi_k, \text{ whenever } i\neq k,\quad i,k=0,\ldots,h\right\},
\end{align*}
see \cite{schnurr:dehling:2017}, Sec. 2.1.\newline\newline
The number of permutations in $S_h$ is given by $\#S_h=(h+1)!$. 
In order to get a better intuitive understanding of the concept of ordinal patterns, we have a closer look at the following example, before turning to the formal definition.
\begin{example}\label{ex:ordinal pattern}
Figure \ref{fig: example ordinal pattern} provides an illustrative understanding of the extraction of an ordinal pattern from a data set. The data points of interest are colored in red and we consider a pattern of length $h=3$, which means we have to take $n=4$ data points into consideration. We fix the points in time $t_0$, $t_1$, $t_2$ and $t_3$ and extract the data points from the time series. Then, we search for the point in time which exhibits the largest value in the resulting data and write down the corresponding time index. In this example it is given by $t=t_1$. We order the data points by writing the time position of the largest value as first entry, the time position of the second largest as second entry and so on. Hence, the absolute values are ordered from largest to smallest and the ordinal pattern $(1,0,3,2)\in S_3$ is obtained for the considered data points.
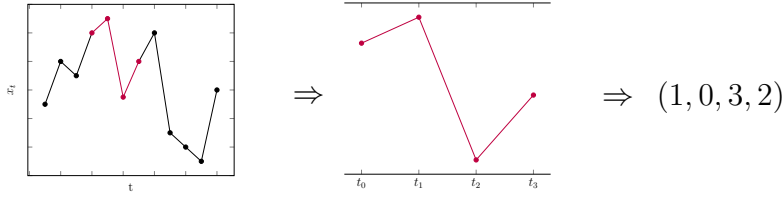
\begin{figure}[htb]
\centering
\begin{minipage}{0.21\textwidth}
\begin{tikzpicture}[scale=0.4]
\begin{axis}[y label style = {at={(axis description cs:0.1,0.5)}},x label style = {at={(axis description cs:0.5,0.075)}},xlabel=t,ylabel=$x_t$,yticklabels={,,},xticklabels={,,}]
\addplot[color=black,mark=*] coordinates{(1,5) (2,8) (3,7) };
\addplot[color=black] coordinates{(3,7)(4,10)};
\addplot[color=black,mark=*] coordinates{(8,10) (9,3) (10,2) (11,1)(12,6)};
\addplot[color=black] coordinates{(7,8) (8,10)};
\addplot[color=purple,mark=*] coordinates{(4,10) (5,11) (6,5.5) (7,8) };
\end{axis}
\end{tikzpicture}
\end{minipage}
\begin{minipage}{0.1\textwidth}
\quad\quad $\Rightarrow$
\end{minipage}
\begin{minipage}{0.21\textwidth}
\begin{tikzpicture}[scale=0.4]
\begin{axis}[xtick={4,5,6,7},xticklabels={$t_0$,$t_1$,$t_2$,$t_3$},hide y axis]
\addplot[color=purple,mark=*] coordinates{(4,10) (5,11) (6,5.5) (7,8) };
\node[
coordinate,
pin = {[rotate=0]$t_0$}
] at (axis cs:0,0) { };
\end{axis}
\end{tikzpicture}
\end{minipage}
\begin{minipage}{0.21\textwidth}
\quad $\Rightarrow\:\:(1,0,3,2)$
\end{minipage}
\caption{Example of the extraction of an ordinal pattern of a given data set.}
\label{fig: example ordinal pattern}
\end{figure}
\end{example}
Formally, the aforementioned procedure can be defined as follows, see \cite{schnurr:dehling:2017}, Section 2.1.
\begin{definition}
As ordinal pattern of a vector $x=\left(x_0,\ldots,x_h\right)\in\mathbb{R}^{h+1}$, we define the unique permutation $\pi=\left(\pi_0,\ldots,\pi_h\right)\in S_h$, 
\begin{align*}
\Pi(x)=\Pi\left(x_0,\ldots,x_h\right)=\left(\pi_0,\ldots,\pi_h\right),
\end{align*}
such that
\begin{align*}
x_{\pi_0}\geq \ldots \geq x_{\pi_h},
\end{align*}
with $\pi_{i-1}<\pi_i$ if $x_{\pi_{i-1}}=x_{\pi_i}$, $i=1,\ldots,h$.
\end{definition}
The last condition assures the uniqueness of $\pi$ if there are ties in the data sets. In particular, this condition is necessary to be made if real-world data is considered.
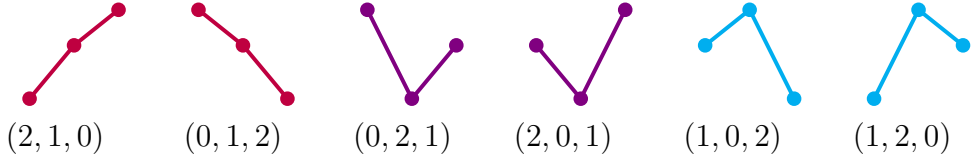
\begin{figure}[h]
\centering
\begin{tikzpicture}
[scale =1]
\begin{axis}[ height =3cm,
  width=3cm,
axis x line=none,
axis y line=none
]
\addplot [mark=*,purple,line width =1.5pt] coordinates {
(0,0.1)(1,0.4)(2,0.6)};
\end{axis}
\end{tikzpicture}\quad\quad
\begin{tikzpicture}
[scale =1]
\begin{axis}[ height =3cm,
  width=3cm,
axis x line=none,
axis y line=none
]
\addplot [mark=*,purple,line width =1.5pt] coordinates {
(0,0.6)(1,0.4)(2,0.1)};
\end{axis}
\end{tikzpicture}\quad\quad
\begin{tikzpicture}
[scale =1]
\begin{axis}[ height =3cm,
  width=3cm,
axis x line=none,
axis y line=none
]
\addplot [mark=*,violet,line width =1.5pt] coordinates {
(0,0.6)(1,0.1)(2,0.4)};
\end{axis}
\end{tikzpicture}\quad\quad
\begin{tikzpicture}
[scale =1]
\begin{axis}[ height =3cm,
  width=3cm,
axis x line=none,
axis y line=none
]
\addplot [mark=*,violet,line width =1.5pt] coordinates {
(0,0.4)(1,0.1)(2,0.6)};
\end{axis}
\end{tikzpicture}\quad\quad
\begin{tikzpicture}
[scale =1]
\begin{axis}[ height =3cm,
  width=3cm,
axis x line=none,
axis y line=none
]
\addplot [mark=*,cyan,line width =1.5pt] coordinates {
(0,0.4)(1,0.6)(2,0.1)};
\end{axis}
\end{tikzpicture}\quad\quad
\begin{tikzpicture}
[scale =1]
\begin{axis}[ height =3cm,
  width=3cm,
axis x line=none,
axis y line=none
]
\addplot [mark=*,cyan,line width =1.5pt] coordinates {
(0,0.1)(1,0.6)(2,0.4)};
\end{axis}
\end{tikzpicture}\newline
$(2,1,0)$\quad\quad\thickspace\thickspace$(0,1,2)$\quad\quad\thickspace$(0,2,1)$\quad\quad$(2,0,1)$\quad\quad\thickspace$(1,0,2)$\quad\thickspace\quad$(1,2,0)$\quad\quad\quad
\caption{Ordinal patterns for $h=2$.}
\label{fig: OP2}
\end{figure}\quad\newline
\noindent In Figure \ref{fig: OP2}, all ordinal patterns of length $h=2$ are shown. As already mentioned in the introduction, from the practical point of view, a highly desirable property of ordinal patterns is that they are not affected by monotone transformations, see \cite{sinn:keller:2011}, p. 1783. \newline
Mathematically, this means: if $f:\mathbb{R}\rightarrow\mathbb{R}$ is strictly montone, then 
\begin{align}
\Pi\left(x_0,\ldots,x_h\right)=\Pi\left(f\left(x_0\right),\ldots,f\left(x_h\right)\right).\label{eq:lineartrafoOP}
\end{align}
In particular, this includes linear transformations $f(x)=ax+b$, with $a\in\mathbb{R}^+$ and $b\in\mathbb{R}$.\newline\newline
Following \cite{schnurr:dehling:2017}, Sec. 1, the minimal requirement of the data sets we use for ordinal analysis in the time series context, i.e., for ordinal pattern probabilities as well as for ordinal pattern dependence later on, is \textit{ordinal pattern stationarity (of order $h$)}. This property implies that the probability of observing a certain ordinal pattern of length $h$ remains the same when shifting the moving window of length $h$ through the entire time series and is not depending on the specific points in time. In the course of this work, the time series, in which the ordinal patterns occur, always have either stationary increments or are even stationary themselves. Note that both properties imply ordinal pattern stationarity. The reason why requiring stationary increments is a sufficient condition is given in the following explanation.\newline\newline
One fundamental property of ordinal patterns is that they are uniquely determined by the increments of the considered time series. As one can imagine in Example \ref{ex:ordinal pattern}, the knowledge of the increments between the data points is sufficient to obtain the corresponding ordinal pattern. In mathematical terms, we can define another mapping $\tilde{\Pi}$, that assigns the corresponding ordinal pattern to each vector of increments, see \cite{sinn:keller:2011}, p. 1783.
\begin{definition}\label{def:tildepi}
We define for $y=\left(y_1,\ldots,y_h\right)\in\mathbb{R}^h$ the mapping $\tilde{\Pi}:\mathbb{R}^h\rightarrow S_h$,
\begin{align*}
\tilde{\Pi}\left(y_1,\ldots,y_h\right):=\Pi\left(0,y_1,y_1+y_2,\ldots,y_1+\ldots+y_h\right),
\end{align*}
such that for $y_i=x_i-x_{i-1}$, $i=1,\ldots,h$, we obtain
\begin{align*}
\tilde{\Pi}\left(y_1,\ldots,y_h\right)&=\Pi\left(0,y_1,y_1+y_2,\ldots,y_1+\ldots+y_h\right)\\
&=\Pi\left(0,x_1-x_0,x_2-x_0,\ldots,x_h-x_0\right)\\
&=\Pi\left(x_0,x_1,x_2,\ldots,x_h\right).
\end{align*}
\end{definition}
We define the two mappings, following \cite{sinn:keller:2011}, p. 1784:
\begin{align*}
&\mathcal{S}:S_h\rightarrow S_h, \left(\pi_0,\ldots,\pi_h\right)\rightarrow\left(\pi_h,\ldots,\pi_0\right),\\
&\mathcal{T}:S_h\rightarrow S_h, \left(\pi_0,\ldots,\pi_h\right)\rightarrow\left(h-\pi_0,\ldots,h-\pi_h\right).
\end{align*}
\begin{figure}[h]
\centering
\begin{tikzpicture}
[scale =1]
\begin{axis}[ height =3cm,
  width=10cm,
axis x line=none,
axis y line=none
]
\addplot [mark=*,purple,line width =1.5pt] coordinates {
(2,0.1)(3,0.6)(4,0.4)(5,0.5)};
\addplot [mark=*,purple,line width =1.5pt] coordinates {
(7,0.5)(8,0.4)(9,0.6)(10,0.1)};
\end{axis}
\end{tikzpicture}
\newline
 $\pi=(1,3,2,0)$  \quad \quad \quad $\mathcal{T}(\pi)=(2,0,1,3)$ \quad \quad  \quad \quad \quad\quad
\end{figure}
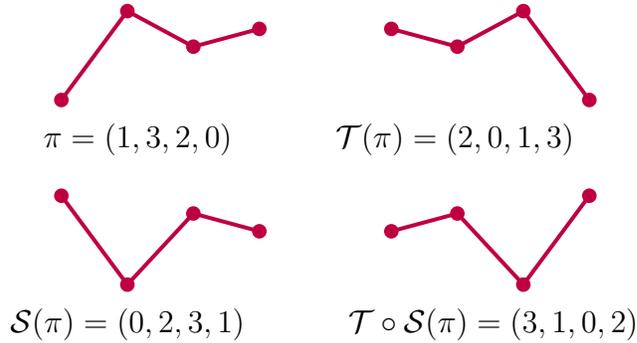
\begin{figure}[h]
\centering
\begin{tikzpicture}
[scale =1]
\begin{axis}[ height =3cm,
  width=10cm,
axis x line=none,
axis y line=none
]
\addplot [mark=*,purple,line width =1.5pt] coordinates {
(2,0.6)(3,0.1)(4,0.5)(5,0.4)};
\addplot [mark=*,purple,line width =1.5pt] coordinates {
(7,0.4)(8,0.5)(9,0.1)(10,0.6)};
\end{axis}
\end{tikzpicture}
\newline
$\mathcal{S}(\pi)=(0,2,3,1)$ \quad \quad \quad   $\mathcal{T}\circ\mathcal{S}(\pi)=(3,1,0,2)$  \quad \quad \quad  \quad \quad 
\caption{Space and time reversion of the pattern $\pi=(1,3,2,0)$.}
\label{fig: setpibar}
\end{figure}~\newline\newline
\noindent
An illustrative understanding of these mappings is given as follows. The mapping $\mathcal{S}(\pi)$, which is the spatial reversion of the pattern $\pi$, is the reflection of $\pi$ on a horizontal line, while $\mathcal{T}(\pi)$, the time reversal of $\pi$, is its reflection on a vertical line, as one can observe in Figure \ref{fig: setpibar}.\newline
Based on the spatial reversion, we define a possibility to divide $S_h$ into two disjoint sets.
\begin{definition}\label{def: Shstern}
We define $S^*_h$ as a subset of $S_h$ with the property that for each $\pi\in S_h$ either $\pi$ or $\mathcal{S}(\pi)$ are contained in the set, but not both of them. 
\end{definition}
Note that this definition does not yield uniqueness of $S_h^*$.\newline\newline
\begin{example}\label{bsp:Sstern} We consider the case $h=2$ again and we want to divide $S_2$ into a possible choice of $S_2^*$ and the corresponding spatial reversal. We choose $S_2^*=\{(2,1,0),(2,0,1),(1,2,0)\}$ and, therefore, $S_2\setminus S_2^*=\{(0,1,2),(1,0,2),(0,2,1)\}$. Remark that $S_2^*=\{(0,1,2),(2,0,1),(1,2,0)\}$ is also a possible choice. The only condition that has to be satisfied is that if one permutation is chosen for $S_2^*$, then its spatial reverse must not be an element of this set.
\end{example}
We stick to the formal definition of ordinal pattern dependence, as it is proposed in \cite{schnurr:dehling:2017}, Sec. 2.1.
The considered moving window consists of $h+1$ data points and, hence, $h$ increments. We define
\begin{align}
p:&=p_{X^{(1)},X^{(2)}}:=\mathbb{P}\left(\Pi\left(X_0^{(1)},\ldots, X_h^{(1)}\right)=\Pi\left(X_0^{(2)},\ldots, X_h^{(2)}\right)\right)\label{eq:p}
\intertext{and}
q:&=q_{X^{(1)},X^{(2)}}:=\sum_{\pi\in S_h} \mathbb{P}\left(\Pi\left(X_0^{(1)},\ldots, X_h^{(1)}\right)=\pi\right)\mathbb{P}\left(\Pi\left(X_0^{(2)},\ldots, X_h^{(2)}\right)=\pi\right).\notag
\end{align}
Then, we define ordinal pattern dependence $OPD$ as
\begin{align}
OPD:=OPD_{X^{(1)},X^{(2)}}:=\frac{p-q}{1-q}.\label{def: OPD}
\end{align}
The parameter $q$ represents the hypothetical case of independence between the two time series. In this case $p$ and $q$ would obtain equal values and, therefore, $OPD$ would equal zero. Regarding the other extreme, the case in which both processes coincide or one is a strictly monotone increasing transform of the other one, we obtain the value $1$. 
However, in the following, we assume $p\in(0,1)$ and $q\in (0,1)$.\newline
Note that the definition of ordinal pattern dependence in \eqref{def: OPD} only measures positive dependence. This is no restriction in practice, because negative dependence can be investigated in an analogous way, by considering $OPD_{X^{(1)},-X^{(2)}}$. If one is interested in both types of dependence simultaneously, in \cite{schnurr:dehling:2017} the authors propose to use $\left(OPD_{X^{(1)},X^{(2)}}\right)_+-\left(OPD_{X^{(1)},-X^{(2)}}\right)_+$. To keep the notation simple, we focus on $OPD$ as it is defined in \eqref{def: OPD}.\newline
\begin{figure}[htb]
\begin{minipage}{0.21\textwidth}
\begin{tikzpicture}[scale=0.35]
\begin{axis}[y label style = {at={(axis description cs:0.1,0.5)}} ,ylabel=$X_t^{(1)}$,yticklabels={,,},xticklabels={,$t_0$,$t_2$,$t_4$,$t_6$,$\ldots$,$t_{n-2}$,$t_{n}$}]
\addplot[color=purple,mark=*] coordinates{(0,3)(1,5) (2,8)  };
\addplot[color=black, no marks] coordinates{(2,8)(3,7)};
\addplot[color=black,mark=*] coordinates{(3,7)(4,10) (5,11) (6,4.5) (7,3) };
\addplot[color=black, no marks] coordinates{(7,3)(8,2.66)};
\addplot[color=black,dotted,thick] coordinates{(8,2.66)(9,2.33) };
\addplot[color=black, no marks] coordinates{(9,2.33)(10,2)};
\addplot[color=black,mark=*] coordinates{(10,2) (11,1)(12,6)};
\end{axis}
\end{tikzpicture}
\end{minipage}
\begin{minipage}{0.21\textwidth}
\begin{tikzpicture}[scale=0.35]
\begin{axis}[y label style = {at={(axis description cs:0.1,0.5)}}, ylabel=$X_t^{(1)}$,yticklabels={,,},xticklabels={,$t_0$,$t_2$,$t_4$,$t_6$,$\ldots$,$t_{n-2}$,$t_{n}$}]
\addplot[color=black,mark=*] coordinates{(0,3)};
\addplot[color=black] coordinates{(0,3)(1,5)};
\addplot[color=purple,mark=*] coordinates{(1,5) (2,8) (3,7)  };
\addplot[color=black] coordinates{(3,7)(4,10)};
\addplot[color=black,mark=*] coordinates{(4,10) (5,11) (6,4.5) (7,3) };
\addplot[color=black, no marks] coordinates{(7,3)(8,2.66)};
\addplot[color=black,dotted,thick] coordinates{(8,2.66)(9,2.33) };
\addplot[color=black, no marks] coordinates{(9,2.33)(10,2)};
\addplot[color=black,mark=*] coordinates{(10,2) (11,1)(12,6)};
\end{axis}
\end{tikzpicture}
\end{minipage}
\begin{minipage}{0.2\textwidth}
\begin{tikzpicture}[scale=0.35]
\begin{axis}[y label style = {at={(axis description cs:0.1,0.5)}} ,ylabel=$X_t^{(1)}$,yticklabels={,,},xticklabels={,$t_0$,$t_2$,$t_4$,$t_6$,$\ldots$,$t_{n-2}$,$t_{n}$}]
\addplot[color=black,mark=*] coordinates{(0,3)(1,5)};
\addplot[color=black] coordinates{(1,5)(2,8)};
\addplot[color=purple,mark=*] coordinates{(2,8) (3,7)(4,10) };
\addplot[color=black] coordinates{(4,10) (5,11)};
\addplot[color=black,mark=*] coordinates{(5,11) (6,4.5) (7,3) };
\addplot[color=black, no marks] coordinates{(7,3)(8,2.66)};
\addplot[color=black,dotted,thick] coordinates{(8,2.66)(9,2.33) };
\addplot[color=black, no marks] coordinates{(9,2.33)(10,2)};
\addplot[color=black,mark=*] coordinates{(10,2) (11,1)(12,6)};\end{axis}
\end{tikzpicture}
\end{minipage}
\begin{minipage}{0.15\textwidth}
\quad\quad\quad\ldots
\end{minipage}
\begin{minipage}{0.2\textwidth}
\begin{tikzpicture}[scale=0.35]
\begin{axis}[y label style = {at={(axis description cs:0.1,0.5)}} ,ylabel=$X_t^{(1)}$,yticklabels={,,},xticklabels={,$t_0$,$t_2$,$t_4$,$t_6$,$\ldots$,$t_{n-2}$,$t_{n}$}]
\addplot[color=black,mark=*] coordinates{(0,3)(1,5)(2,8)(3,7)(4,10)(5,11) (6,4.5) (7,3) };
\addplot[color=black, no marks] coordinates{(7,3)(8,2.66)};
\addplot[color=black,dotted,thick] coordinates{(8,2.66)(9,2.33) };
\addplot[color=black, no marks] coordinates{(9,2.33)(10,2)};
\addplot[color=purple,mark=*] coordinates{(10,2) (11,1)(12,6)};\end{axis}
\end{tikzpicture}
\end{minipage}
\vspace{0.4cm}
\newline
\begin{minipage}{0.21\textwidth}
\begin{tikzpicture}[scale=0.35]
\begin{axis}[y label style = {at={(axis description cs:0.1,0.5)}} ,ylabel=$X_t^{(2)}$,yticklabels={,,},xticklabels={,$t_0$,$t_2$,$t_4$,$t_6$,$\ldots$,$t_{n-2}$,$t_{n}$}]
\addplot[color=purple,mark=*] coordinates{(0,3)(1,7) (2,8) };
\addplot[color=black] coordinates{(2,8)(3,5)};
\addplot[color=black,mark=*] coordinates{(3,5)(4,10) (5,1) (6,5.5) (7,8) };
\addplot[color=black] coordinates{(7,8) (8,7.33)};
\addplot[color=black,dotted,thick] coordinates{(8,7.33)(9,6.66) };
\addplot[color=black] coordinates{(9,6.66) (10,6)};
\addplot[color=black,mark=*] coordinates{ (10,6) (11,3)(12,1)};
\end{axis}
\end{tikzpicture}
\end{minipage}
\begin{minipage}{0.21\textwidth}
\begin{tikzpicture}[scale=0.35]
\begin{axis}[y label style = {at={(axis description cs:0.1,0.5)}} ,ylabel=$X_t^{(2)}$,yticklabels={,,},xticklabels={,$t_0$,$t_2$,$t_4$,$t_6$,$\ldots$,$t_{n-2}$,$t_{n}$}]
\addplot[color=black,mark=*] coordinates{(0,3)};
\addplot[color=black] coordinates{(0,3)(1,7)};
\addplot[color=purple,mark=*] coordinates{(1,7) (2,8)(3,5) };
\addplot[color=black] coordinates{(3,5)(4,10)};
\addplot[color=black,mark=*] coordinates{(4,10) (5,1) (6,5.5) (7,8) };
\addplot[color=black] coordinates{(7,8) (8,7.33)};
\addplot[color=black,dotted,thick] coordinates{(8,7.33)(9,6.66) };
\addplot[color=black] coordinates{(9,6.66) (10,6)};
\addplot[color=black,mark=*] coordinates{ (10,6) (11,3)(12,1)};
\end{axis}
\end{tikzpicture}
\end{minipage}
\begin{minipage}{0.2\textwidth}
\begin{tikzpicture}[scale=0.35]
\begin{axis}[y label style = {at={(axis description cs:0.1,0.5)}},ylabel=$X_t^{(2)}$,yticklabels={,,},xticklabels={,$t_0$,$t_2$,$t_4$,$t_6$,$\ldots$,$t_{n-2}$,$t_{n}$}]
\addplot[color=black,mark=*] coordinates{(0,3)(1,7)};
\addplot[color=black] coordinates{(1,7)(2,8)};
\addplot[color=purple,mark=*] coordinates{(2,8)(3,5)(4,10) };
\addplot[color=black] coordinates{(4,10)(5,1)};
\addplot[color=black,mark=*] coordinates{ (5,1) (6,5.5) (7,8) };
\addplot[color=black] coordinates{(7,8) (8,7.33)};
\addplot[color=black,dotted,thick] coordinates{(8,7.33)(9,6.66) };
\addplot[color=black] coordinates{(9,6.66) (10,6)};
\addplot[color=black,mark=*] coordinates{ (10,6) (11,3)(12,1)};
\end{axis}
\end{tikzpicture}
\end{minipage}
\begin{minipage}{0.15\textwidth}
\quad\quad\quad\ldots
\end{minipage}
\begin{minipage}{0.2\textwidth}
\begin{tikzpicture}[scale=0.35]
\begin{axis}[y label style = {at={(axis description cs:0.1,0.5)}} ,ylabel=$X_t^{(2)}$,yticklabels={,,},xticklabels={,$t_0$,$t_2$,$t_4$,$t_6$,$\ldots$,$t_{n-2}$,$t_{n}$}]
\addplot[color=black,mark=*] coordinates{(0,3)(1,7)};
\addplot[color=black] coordinates{(1,7)(2,8)};
\addplot[color=black,mark=*] coordinates{(2,8)(3,5)(4,10) };
\addplot[color=black] coordinates{(4,10)(5,1)};
\addplot[color=black,mark=*] coordinates{ (5,1) (6,5.5) (7,8) };
\addplot[color=black] coordinates{(7,8) (8,7.33)};
\addplot[color=black,dotted,thick] coordinates{(8,7.33)(9,6.66) };
\addplot[color=black] coordinates{(9,6.66) (10,6)};
\addplot[color=purple,mark=*] coordinates{ (10,6) (11,3)(12,1)};
\end{axis}
\end{tikzpicture}
\end{minipage}
\caption{Illustration of estimation of ordinal pattern dependence.}
\label{fig:OPD}
\end{figure}
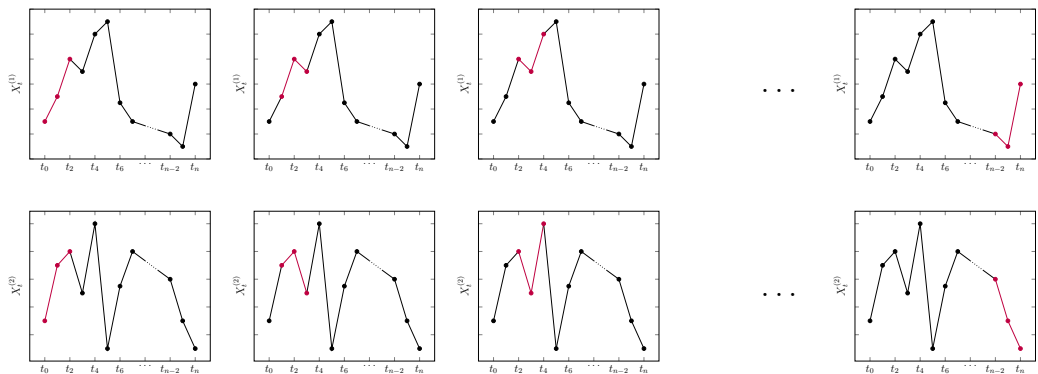\quad\newline\newline
We compare whether the ordinal patterns in $\left(X_j^{(1)}\right)_{j\in\mathbb{Z}}$ coincide with the ones in $\left(X_j^{(2)}\right)_{j\in\mathbb{Z}}$. Recall that it is an essential property of ordinal patterns that they are uniquely determined by the increment process. Therefore, we have to consider the increment processes $\left(Y_j\right)_{j\in\mathbb{Z}}=\left(\left(Y_j^{(1)},Y_j^{(2)}\right)\right)_{j\in\mathbb{Z}}$ as defined in \eqref{eq:process} for $d=2$, where $Y_j^{(p)}=X_j^{(p)}-X_{j-1}^{(p)}$, $p=1,2$. Hence, we can also express $p$ and $q$ (and consequently $OPD$) as a probability that only depends on the increments of the considered vectors of the time series. Recall the definition of $\left(Y_{j,h}\right)_{j\in\mathbb{Z}}$ for $d=2$, given by
\begin{align*}
Y_{j,h}=\left(Y_j^{(1)},\ldots,Y_{j+h-1}^{(1)},Y_j^{(2)},\ldots,Y_{j+h-1}^{(2)}\right)^t,
\end{align*}
such that $Y_{j,h}\sim\mathcal{N}\left(0,\Sigma_{2,h}\right)$ with $\Sigma_{2,h}$ as given in \eqref{eq: covariance matrix extended process}.\newline\newline
In the course of this article, we focus on the estimation of $p$. For a detailed investigation of the limit theorems for estimators of $OPD$ we refer to \cite{nuessgen:2021}.
We define the estimator of $p$, the probability of coincident patterns in both time series in a moving window of fixed length, by
\begin{align*}
\hat{p}_n&=\frac{1}{n-h}\sum_{j=0}^{n-h-1}\mathbf{1}_{\left\{\Pi\left(X_{j}^{(1)},\ldots, X_{j+h}^{(1)}\right)=\Pi\left(X_j^{(2)},\ldots, X_{j+h}^{(2)}\right)\right\}}\\
&=\frac{1}{n-h}\sum_{j=1}^{n-h}\mathbf{1}_{\left\{\tilde{\Pi}\left(Y_{j}^{(1)},\ldots, Y_{j+h-1}^{(1)}\right)=\tilde{\Pi}\left(Y_{j}^{(2)},\ldots, Y_{j+h-1}^{(2)}\right)\right\}},
\end{align*}
where 
\begin{align*}
\tilde{\Pi}\left(Y_1,\ldots,Y_h\right):&=\Pi\left(0,Y_1,Y_1+Y_2,\ldots,Y_1+\ldots+Y_h\right)\\
&=\Pi\left(0,X_1-X_0,\ldots,X_h-X_0\right)\\
&=\Pi\left(X_0,X_1,\ldots,X_h\right).
\end{align*}
Figure \ref{fig:OPD} illustrates the way ordinal pattern dependence is estimated by $\hat{p}_n$. The patterns of interest that are compared in each moving window are colored in red.\newline\newline
Having emphasized the crucial importance of the increments, we define the following conditions on the increment process $\left(Y_j\right)_{j \in\mathbb{Z}}$: let $\left(Y_j\right)_{j \in\mathbb{Z}}$ be a bivariate, stationary Gaussian process with $Y_j^{(p)}\sim\mathcal{N}(0,1)$, $p=1,2$.
\begin{enumerate}[label=\textbf{(L\arabic*})]
\renewcommand{\labelenumi}{\textbf{(L)}}
\item
We assume that $\left(Y_j\right)_{j \in\mathbb{Z}}$ fulfills \eqref{multivariateLRDconditionII} with $d^*$ in $\left(\frac{1}{4},\frac{1}{2}\right)$. We allow for $\min\left\{d_1,d_2\right\}$ to be in the range $\left(-\infty,0\right)\cup\left(0,\frac{1}{4}\right]$. \label{it:L}
\end{enumerate}
\begin{enumerate}[label=\textbf{(S})]
\renewcommand{\labelenumi}{\textbf{(S)}}
\item We assume $d_1,d_2\in\left(-\infty,0\right)\cup\left(0,\frac{1}{4}\right)$ such that the cross-correlation function of $\left(Y_j\right)_{j \in\mathbb{Z}}$ fulfills for $p,q=1,2$
\begin{align*}
r^{(p,q)}(k)= k^{d_p+d_q-1}L_{p,q}(k)\quad(k\rightarrow\infty)
\end{align*}
with $L_{p,q}(k)\rightarrow L_{p,q}$ and $L_{p,q}\in\mathbb{R}$ holds.\label{it:S}
\end{enumerate}
Furthermore, in both cases, it holds that $\left|r^{(p,q)}(k)\right|<1$ for $p,q=1,2$ and $k\geq 1$ to exclude ties.\newline\newline
We begin with the investigation of the asymptotics of $\hat{p}_n$. First, we calculate the Hermite rank of $\hat{p}_n$, since the Hermite rank determines for which ranges of $d^*$ the estimator $\hat{p}_n$ is still long-range dependent. Depending on this range, different limit theorems may hold.
\begin{Lem}\label{lemm:OPDHR2}
The Hermite rank of $f(Y_{j,h})=\mathbf{1}_{\left\{\tilde{\Pi}\left(Y_{j+1}^{(1)},\ldots, Y_{j+h}^{(1)}\right)=\tilde{\Pi}\left(Y_{j+1}^{(2)},\ldots, Y_{j+h}^{(2)}\right)\right\}}$ with respect to $\Sigma_{2,h}$ is equal to $2$.
\end{Lem}
\begin{pf}
Following \cite{betken:2019}, Lemma 5.4 it is sufficient to show the following two properties:
\renewcommand{\labelenumi}{(\roman{enumi})}
\begin{enumerate}
\item $m(f,\Sigma_{2,h})\geq 2$,
\item $m(f,I_{2,h})\leq 2$.
\end{enumerate}
Note that the conclusion is not trivial, because $m(f,\Sigma_{2,h})\neq m(f,I_{2,h})$ in general, see \cite{beran:feng:ghosh:kulik:2013}, Lemma 3.7. 
Lemma 5.4 in \cite{betken:2019} can be applied due to the following reasoning. Ordinal patterns are not affected by scaling, therefore, the technical condition that $\Sigma_{2,h}^{-1}-I_{2,h}$ is positive semidefinite is fulfilled in our case. We can scale the standard deviation of the random vector $Y_{j,h}$ by any positive real number $\sigma>0$ since for all $j\in\mathbb{Z}$ we have 
\begin{align*}
{\left\{\tilde{\Pi}\left(Y_{j}^{(1)},\ldots, Y_{j+h-1}^{(1)}\right)=\tilde{\Pi}\left(Y_j^{(2)},\ldots, Y_{j+h-1}^{(2)}\right)\right\}}\\
={\left\{\tilde{\Pi}\left(\sigma Y_{j}^{(1)},\ldots, \sigma Y_{j+h-1}^{(1)}\right)=\tilde{\Pi}\left(\sigma Y_j^{(2)},\ldots, \sigma Y_{j+h-1}^{(2)}\right)\right\}}.
\end{align*}
To show property $(i)$, we need to consider a multivariate random vector 
\begin{align*}
Y_{1,h}:=\left(Y_1^{(1)},\ldots,Y_h^{(1)},Y_1^{(2)},\ldots,Y_h^{(2)}\right)^t
\end{align*}
with covariance matrix $\Sigma_{2,h}$. We fix $i=1,\ldots,2h$.
We divide the set $S_h$ into disjoint sets, namely into $S_h^*$, as defined in Definition \ref{def: Shstern} and the complimentary set $S_h\setminus S_h^*$. Note that
\begin{align*}
-Y_{j,h}\overset{\mathcal{D}}{=} Y_{j,h}
\end{align*}
holds. This implies
\begin{align*}
\mathbb{E}\left(Y_{j,h}^{(i)}\mathbf{1}_{\left\{\tilde{\Pi}\left(Y_{1}^{(1)},\ldots, Y_{h}^{(1)}\right)=\tilde{\Pi}\left(Y_{1}^{(2)},\ldots, Y_{h}^{(2)}\right)=\pi\right\}}\right)=-\mathbb{E}\left(Y_{j,h}^{(i)}\mathbf{1}_{\left\{\tilde{\Pi}\left(Y_{1}^{(1)},\ldots, Y_{h}^{(1)}\right)=\tilde{\Pi}\left(Y_{1}^{(2)},\ldots, Y_{h}^{(2)}\right)=\mathcal{S}(\pi)\right\}}\right)
\end{align*} for $\pi\in S_h$. Hence, we arrive at:
\begin{align*}
\mathbb{E}\left(Y_{j,h}^{(i)}f(Y_{j,h})\right)&=\mathbb{E}\left(Y_{j,h}^{(i)}\mathbf{1}_{\left\{\tilde{\Pi}\left(Y_{1}^{(1)},\ldots, Y_{h}^{(1)}\right)=\tilde{\Pi}\left(Y_{1}^{(2)},\ldots, Y_{h}^{(2)}\right)\right\}}\right)\\
&=\sum_{\pi\in S_h}\mathbb{E}\left(Y_{j,h}^{(i)}\mathbf{1}_{\left\{\tilde{\Pi}\left(Y_{1}^{(1)},\ldots, Y_{h}^{(1)}\right)=\tilde{\Pi}\left(Y_{1}^{(2)},\ldots, Y_{h}^{(2)}\right)=\pi\right\}}\right)\\
&=\sum_{\pi\in S_h^*}\mathbb{E}\left(Y_{j,h}^{(i)}\mathbf{1}_{\left\{\tilde{\Pi}\left(Y_{1}^{(1)},\ldots, Y_{h}^{(1)}\right)=\tilde{\Pi}\left(Y_{1}^{(2)},\ldots, Y_{h}^{(2)}\right)=\pi\right\}}\right)
\\&\quad\quad-\sum_{\pi\in S_h\setminus S_h^*}\mathbb{E}\left(Y_{j,h}^{(i)}\mathbf{1}_{\left\{\tilde{\Pi}\left(Y_{1}^{(1)},\ldots, Y_{h}^{(1)}\right)=\tilde{\Pi}\left(Y_{1}^{(2)},\ldots, Y_{h}^{(2)}\right)=\mathcal{S}(\pi)\right\}}\right)\\
&=0
\end{align*}
for $i=1,\ldots,2h$.\newline
Consequently $m\left(f,\Sigma_{2,h}\right)\geq 2$.\newline\newline
In order to proof $(ii)$, we consider
\begin{align*}
U_{1,h}:=\left(U_1^{(1)},\ldots,U_h^{(1)},U_1^{(2)},\ldots,U_h^{(2)}\right)^t
\end{align*}
 to be a random vector with independent $\mathcal{N}(0,1)$ distributed entries. For $i=1,\ldots,h$ and $k=h+1,\ldots,2h$ such that $k-h=i$, we obtain
\begin{align*}
\mathbb{E}\left(U_{1,h}^{(i)}U_{1,h}^{(k)}f\left(U_{1,h}\right)\right)&=\mathbb{E}\left(U_{i}^{(1)}U_{k-h}^{(2)}\mathbf{1}_{\left\{\tilde{\Pi}\left(U_{1}^{(1)},\ldots, U_{h}^{(1)}\right)=\tilde{\Pi}\left(U_{1}^{(2)},\ldots, U_{h}^{(2)}\right)\right\}}\right)\\
&=\sum_{\pi\in S_h} \mathbb{E}\left(U_{i}^{(1)}U_{i}^{(2)}\mathbf{1}_{\left\{\tilde{\Pi}\left(U_{1}^{(1)},\ldots, U_{h}^{(1)}\right)=\tilde{\Pi}\left(U_{1}^{(2)},\ldots, U_{h}^{(2)}\right)=\pi\right\}}\right)\\
&=\sum_{\pi\in S_h} \left( \mathbb{E}\left(  U_i^{(1)} \mathbf{1}_{\left\{  \tilde{\Pi}\left(U_1^{(1)},\ldots,U_{h}^{(1)}\right)=\pi  \right\}}   \right)\right)^2\\
&\neq 0,
\end{align*}
since $\mathbb{E}\left(U_{i}^{(1)}\mathbf{1}_{\left\{  \tilde{\Pi}\left(U_1^{(1)},\ldots,U_{h}^{(1)}\right)=\pi  \right\}}   \right)\neq 0$ for all $\pi\in S_h$. This was shown in the proof of Lemma 3.4 in \cite{betken:2019}.\newline
All in all, we derive $m(f,\Sigma_{2,h})=2$ and, hence, have proven the lemma. $\hfill\Box$\newline
\end{pf}\quad\newline
The case $m(f,\Sigma_{2,h})=2$ exhibits the property that the standard range of the long-range dependence parameter $d^*\in\left(0,\frac{1}{2}\right)$ has to be divided into two different sets. If $d^*\in\left(\frac{1}{4},\frac{1}{2}\right)$, the transformed process $f\left(Y_{j,h}\right)_{j\in\mathbb{Z}}$ is still long-range dependent, see \cite{pipiras:taqqu:2017}, Table 5.1. If $d^*\in\left(0,\frac{1}{4}\right)$, the transformed process is short-range dependent, which means by definition that the autocorrelations of the transformed process are summable, see \cite{kechagias:pipiras:2015}, Remark 2.3. Therefore, we have two different asymptotic distributions that have to be considered for the estimator $\hat{p}_n$ of coincident patterns.
\subsection{Limit theorem for the estimator of $p$ in case of long-range dependence}\label{sec: LTpLRD}
First, we restrict ourselves to the case that at least one of the two parameters $d_1$ and $d_2$ is in $\left(\frac{1}{4},\frac{1}{2}\right)$. This assures $d^*\in\left(\frac{1}{4},\frac{1}{2}\right)$. We explicitly include mixing cases where the process corresponding to $\min\left\{d_1,d_2\right\}$ is allowed to be long-range as well as short-range dependent. \newline
Note that this setting includes the pure long-range dependence case, which means that for $p=1,2$, we have $d_p\in\left(\frac{1}{4},\frac{1}{2}\right)$, or even $d_1=d_2=d^*$. However, in general the assumptions are lower, such that we only require $d_p\in\left(\frac{1}{4},\frac{1}{2}\right)$ for either $p=1$ or $p=2$ and the other parameter is allowed to be in $\left(-\infty,0\right)$ or $\left(0,\frac{1}{4}\right)$, too.\newline\newline
We can, therefore, apply the results of Corollary \ref{cor:LTHR2} and obtain the following asymptotic distribution for $\hat{p}_n$:
\begin{theorem}\label{LT:OPDp}
Under the assumption in \ref{it:L}, we obtain
\begin{align}
n^{1-2d^*}(C_2)^{-\frac{1}{2}}&\left(\hat{p}_n-p\right)\xrightarrow{\mathcal{D}} \sum_{p,q\in P^*}\tilde{\alpha}^{(p,q)} Z^{(p,q)}_{2,d^*+1/2}(1)
\end{align}
with $Z^{(p,q)}_{2,d^*+1/2}(1)$ as given in Theorem \ref{th: Hermite 2} for $p,q\in P^*$ and $C_2:=\frac{1}{2d^*\left(4d^*-1\right)}$ being a normalizing constant. We have 
\begin{align*}
\tilde{\alpha}^{(p,q)}:=\sum\limits_{i,k=1}^h \alpha_{i,k}^{(p,q)}, \text{   where  } \alpha^{(p,q)}_{i,k}=\alpha_{i+(p-1)h,k+(q-1)h},
\end{align*}
for each $p,q\in P^*$ and $i,k=1,\ldots,h$ and $\left(\alpha_{i,k}\right)_{1\leq i,k\leq dh}=\Sigma_{2,h}^{-1}C\Sigma_{2,h}^{-1}$, where the variable 
\begin{align*}
C=\left(c_{i,k}\right)_{1\leq i,k\leq 2h}=\mathbb{E}\left(Y_{1,h}\left( \mathbf{1}_{\left\{\tilde{\Pi}\left(Y_{1}^{(1)},\ldots, Y_{h}^{(1)}\right)=\tilde{\Pi}\left(Y_1^{(2)},\ldots, Y_{h}^{(2)}\right)\right\}}-p\right) Y_{1,h}^t\right)
\end{align*}
denotes the matrix of second order Hermite coefficients.
\end{theorem}
\begin{pf}
The proof of this theorem is an immediate application of Corollary \ref{cor:LTHR2} and Lemma \ref{lemm:OPDHR2}. Note that for $\hat{p}_n$ it holds that it is square integrable with respect to $Y_{j,h}$ and that the set of discontinuity points is a Null set with respect to the $2h$-dimensional Lebesgue measure. This is shown in \cite{nuessgen:2021}, eq. (4.5).\newline\newline.$\hfill\Box$\newline
\end{pf}\quad\newline
Following Theorem \ref{th: HR2alternative}, we are also able to express the limit distribution above in terms of two standard Rosenblatt random variables by modifying the weighting factors in the limit distribution. Note that this requires slightly stronger assumptions as in Theorem \ref{th: Hermite 2}.
\begin{theorem}\label{LT:OPDpAlternative}
Let \ref{it:L} hold with $d_1=d_2$. Additionally we assume that $r^{(1,1)}(l)=r^{(2,2)}(l)$, for $l=0,\ldots,h-1$, and $L_{1,1}+L_{2,2}\neq L_{1,2}+L_{2,1}$. Then we obtain
\begin{align*}
n^{1-2d^*}(C_2)^{-\frac{1}{2}}\left(\hat{p}_n-p\right)&\xrightarrow{\mathcal{D}} 
 \left(\tilde{\alpha}^{(1,1)}-\tilde{\alpha}^{(1,2)}\right) \frac{L_{2,2}-L_{2,1}-L_{1,2}+L_{1,1}}{2}Z^*_{2,d^*+1/2}(1)\\
&\quad\quad+\left(\tilde{\alpha}^{(1,1)}+\tilde{\alpha}^{(1,2)}\right) \frac{L_{2,2}+L_{2,1}+L_{1,2}+L_{1,1}}{2}Z^{**}_{2,d^*+1/2}(1),
\end{align*}
with $C_2$ and $\tilde{\alpha}^{(p,q)}$ as given in Theorem \ref{LT:OPDp}.
Note that $Z^{*}_{2,d^*+1/2}(1)$ and $Z^{**}_{2,d^*+1/2}(1)$ are both standard Rosenblatt random variables whose covariance is given by
\begin{align}
Cov\left(Z^{*}_{2,d^*+1/2}(1),Z^{**}_{2,d^*+1/2}(1)\right)=\frac{\left(L_{2,2}-L_{1,1}\right)^2  }{\left(L_{1,1}+L_{2,2}\right)^2-\left(L_{1,2}+L_{2,1}\right)^2}.\label{eq:covarianceOPDalternative}
\end{align} 
\end{theorem}
\begin{remark}
Following \cite{nuessgen:2021}, Corollary 3.14, if additionally $r^{(1,1)}(k)=r^{(2,2)}(k)$ and $r^{(1,2)}(k)=r^{(2,1)}(k)$ is fulfilled for all $k\in\mathbb{Z}$, then, the two limit random variables following a standard Rosenblatt distribution in Theorem \ref{LT:OPDpAlternative} are independent. Note that due to the considerations in \cite{veillette:taqqu:2013}, eq. (10), we know that the distribution of the sum of two independent standard Rosenblatt random variables is not standard Rosenblatt. However, this yields a computational benefit, as it is possible to efficiently simulate the standard Rosenblatt distribution, for details, see \cite{veillette:taqqu:2013}.
\end{remark}
We turn to an example that deals with the asymptotic variance of the estimator of $p$ in Theorem \ref{LT:OPDp} in the case $h=1$.
\begin{example}\label{ex:asymptvarianceh=1}
We focus on the case $h=1$ and consider the underlying process $\left(Y_{j,1}\right)_{j\in\mathbb{Z}}=\left(Y_j^{(1)},Y_j^{(2)}\right)_{j\in\mathbb{Z}}$. It is possible to determine the asymptotic variance depending on the correlation $r^{(1,2)}(0)$ between these two increment variables.\newline
We start with the calculation of the second order Hermite coefficients in the case $\pi=(1,0)$. This corresponds to the event $\left\{Y_j^{(1)}\geq 0,Y_j^{(2)}\geq 0\right\}$, which yields
\begin{align*}
c_{1,1}^{\pi,2}=\mathbb{E}\left(  \left(\left(Y_j^{(1)}\right)^2-1\right)\mathbf{1}_{\left\{Y_j^{(1)}\geq 0,Y_j^{(2)}\geq 0\right\}}  \right)
\end{align*}
and
\begin{align*}
c_{1,2}^{\pi,2}=\mathbb{E}\left(  \left(Y_j^{(1)}Y_j^{(2)}\right)\mathbf{1}_{\left\{Y_j^{(1)}\geq 0,Y_j^{(2)}\geq 0\right\}}  \right).
\end{align*}
Due to $r^{(1,2)}(0)=r^{(2,1)}(0)$, we have $\left(Y_j^{(1)},Y_j^{(2)}\right)\overset{\mathcal{D}}{=}\left(Y_j^{(2)},Y_j^{(1)}\right)$ and, therefore, $c_{1,1}^{\pi,2}=c_{2,2}^{\pi,2}$. We identify the second order Hermite coeficients as the ones already calculated in \cite{betken:2019}, Example 3.13, although we are considering two consecutive increments of a univariate Gaussian process there. However, since the corresponding values are only determined by the correlation between the Gaussian variables, we can simply replace the autocorrelation at lag $1$ by the cross-correlation at lag $0$. Hence, we obtain
\begin{align*}
c_{1,1}^{\pi,2}&=\varphi^2(0)r^{(1,2)}(0)\sqrt{1-\left(r^{(1,2)}(0)  \right)^2},\\
c_{1,2}^{\pi,2}&=\varphi^2(0)\sqrt{1-\left(r^{(1,2)}(0)  \right)^2}.
\end{align*}
Recall that the inverse $\Sigma_{2,1}^{-1}=\left(g_{i,j}\right)_{i,j=1,2}$ of the correlation matrix of $\left(Y_j^{(1)},Y_j^{(2)}\right)$ is given by
\begin{align*}
\Sigma_{2,1}^{-1}=\frac{1}{1-\left(r^{(1,2)}(0)  \right)^2}\begin{pmatrix}
1 && -r^{(1,2})(0) \\ -r^{(1,2)}(0) && 1
\end{pmatrix}.
\end{align*}
By using the formula for $\tilde{\alpha}^{(p,q)}$ obtained in \cite{nuessgen:2021}, eq. (4.23), we derive
\begin{align*}
\tilde{\alpha}^{(1,1)}_{\pi,2}&=\alpha^{\pi,2}_{1,1}=\left(g_{1,1}^2+g_{1,2}^2\right)c_{1,1}^{\pi,2} + 2g_{1,1}g_{1,2}c_{1,2}^{\pi,2},\\
\tilde{\alpha}^{(1,2)}_{\pi,2}&=\alpha^{\pi,2}_{1,2}=\left(g_{1,1}^2+g_{1,2}^2\right)c_{1,2}^{\pi,2} + 2g_{1,1}g_{1,2}c_{1,1}^{\pi,2}.
\end{align*}
Plugging the second order Hermite cofficients and the entries of the inverse of the covariance matrix depending on $r^{(1,2)}(0)$ into the formulas, we arrive at
\begin{align*}
\tilde{\alpha}^{(1,1)}_{\pi,2}=\frac{-\varphi^2(0)r^{(1,2)}(0)}{\left( 1-\left(r^{(1,2)}(0)\right)^2  \right)^{1/2}}
\end{align*}
and \begin{align*}
\tilde{\alpha}^{(1,2)}_{\pi,2}=\frac{\varphi^2(0)}{\left(  1-\left(r^{(1,2)}(0)\right)^2 \right)^{1/2}}.
\end{align*}
Therefore, in the case $h=1$, we obtain the following factors in the limit variance in Theorem \ref{LT:OPDp}:
\begin{align*}
\tilde{\alpha}^{(1,1)}&=\tilde{\alpha}^{(2,2)}=\frac{-2\varphi^2(0)r^{(1,2)}(0)}{\left(  1-\left(r^{(1,2)}(0)\right)^2  \right)^{1/2}}\\
\tilde{\alpha}^{(1,2)}&=\tilde{\alpha}^{(2,1)}=\frac{2\varphi^2(0)}{\left(  1-\left(r^{(1,2)}(0)\right)^2 \right)^{1/2}}.
\end{align*}
\end{example}
\begin{remark}
It is not possible to determine the limit variance for $h=2$ analytically, as this includes orthant probabilities of a four-dimensional Gaussian disitribution. Following \cite{abrahamson:1964} no closed formulas are available for these probabilities. However, there are fast algorithms at hand that calculate the limit variance efficiently. It is possible to take advantage of the symmetry properties of the multivariate Gaussian distribution to keep the computational cost of these algorithms low. For details, see \cite{nuessgen:2021}, Sec. 4.3.1. 
\end{remark}
\subsection{Limit theorem for the estimator of $p$ in case of short-range dependence}\label{sec:OPDSRD}
In this section, we focus on the case that $d^*\in \left(-\infty,0\right)\cup\left(0,\frac{1}{4}\right)$. If $d^*\in \left(0,\frac{1}{4}\right)$, we are still dealing with a long-range dependent multivariate Gaussian process $\left(Y_{j,h}\right)_{j\in\mathbb{Z}}$. However, the transformed process $\hat{p}_n-p$ is no longer long-range dependent, since we are considering a function with Hermite rank $2$, see also \cite{pipiras:taqqu:2017}, Table 5.1. Otherwise, if $d^*\in \left(-\infty,0\right)$, the process $\left(Y_{j,h}\right)_{j\in\mathbb{Z}}$ itself is already short-range dependent, since the cross-correlations are summable. Therefore, we obtain the following central limit theorem by applying Theorem 4 in \cite{arcones:1994}.
\begin{theorem}\label{LT: SRDp}
Under the assumptions in \ref{it:S}, we obtain
\begin{align*}
n^{\frac{1}{2}}\left(\hat{p}_n-p\right)\xrightarrow{\mathcal{D}} \mathcal{N}\left(0,\sigma^2\right)
\end{align*}
with 
\begin{align*}
\sigma^2&=\sum_{k=-\infty}^{\infty} \mathbb{E} \Bigg[ \left(\mathbf{1}_{\left\{\tilde{\Pi}\left(Y_{1}^{(1)},\ldots, Y_{h}^{(1)}\right)=\tilde{\Pi}\left(Y_{1}^{(2)},\ldots, Y_{h}^{(2)}\right)\right\}}-p \right)\\ &\quad\quad\quad\quad\quad\quad\times\left(\mathbf{1}_{\left\{\tilde{\Pi}\left(Y_{1+k}^{(1)},\ldots, Y_{h+k}^{(1)}\right)=\tilde{\Pi}\left(Y_{1+k}^{(2)},\ldots, Y_{h+k}^{(2)}\right)\right\}}-p\right)\Bigg].
\end{align*}
\end{theorem}
We close this section with a brief retrospect of the results obtained. We established limit theorems for the estimator of $p$ as probability of coincident pattern in both time series and, hence, on the most important parameter in the context of ordinal pattern dependence. The long-range dependent case as well as the mixed case of short- and long-range dependence was considered. Finally, we provided a central limit theorem for a multivariate Gaussian time series that is short-range dependent if transformed by $\hat{p}_n$. In the subsequent section, we provide a simulation study that illustrate our theoretical findings. In doing so, we shed light on the Rosenblatt distribution and the distribution of the sum of Rosenblatt distributed random variables.
\section{Simulation study}
\begin{figure}[h!]
\centering
\includegraphics[width=0.8\textwidth]{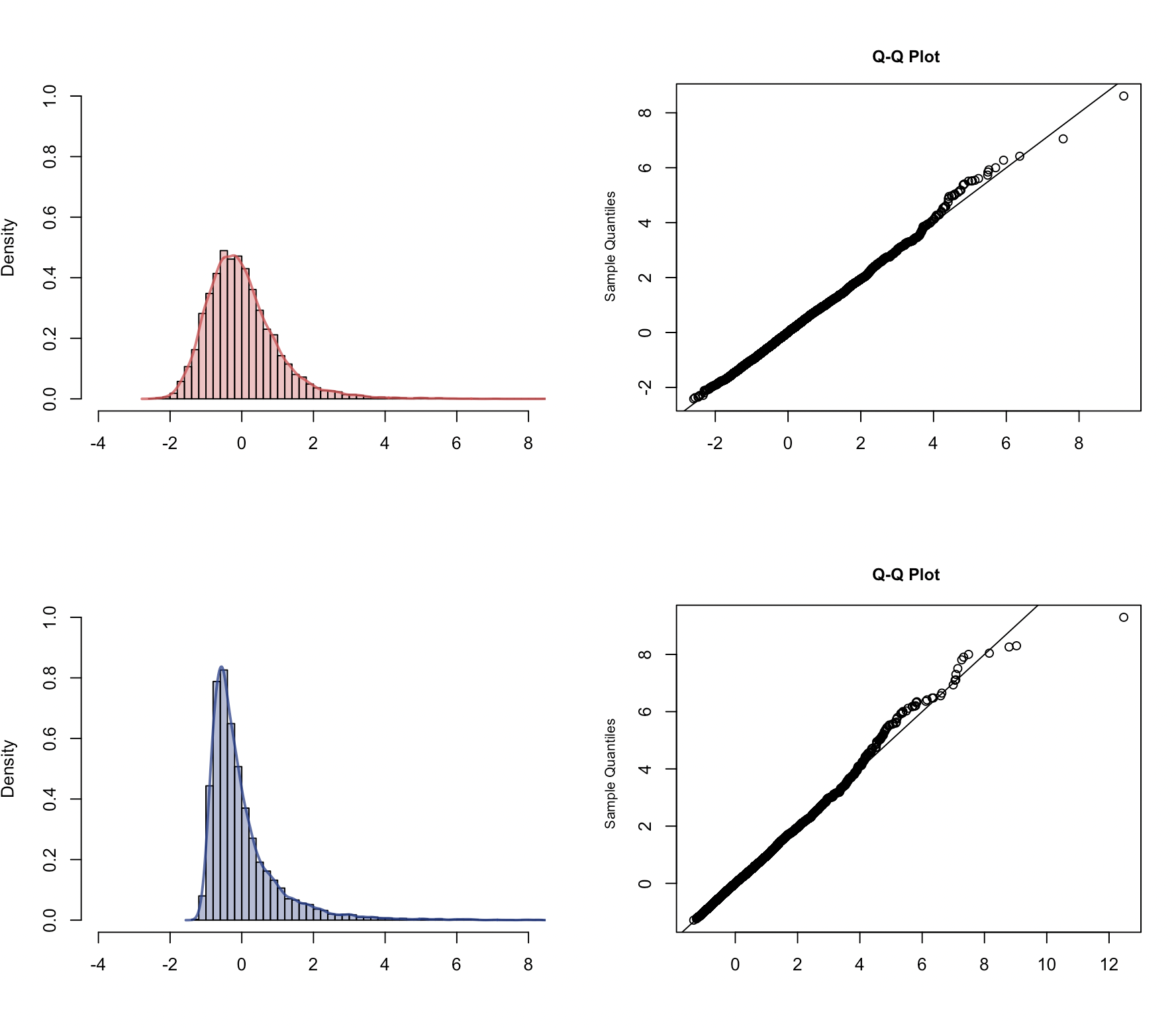}
\caption{Histogram, kernel density estimation and Q-Q plot with respect to the Rosenblatt distribution of $\frac{1}{n}\sum_{j=1}^n H_2\left(Y^*_j\right)$ for different Hurst parameters: $H=0.8$ (top), $H=0.9$ (bottom).}
\label{fig: OPDalt}
\end{figure}

We begin with the generation of a bivariate long-range dependent fractional Gaussian noise series $\left(Y_j^{(1)},Y_j^{(2)}\right)_{j=1,\ldots,n}$. First, we simulate two independent fractional Gaussian noise processes $\left(U_j^{(1)}\right)_{j=1,\ldots,n}$ and $\left(U_j^{(2)}\right)_{j=1,\ldots,n}$ derived by the R-package ``longmemo'', for a fixed parameter $H\in\left(\frac{1}{2},1\right)$ in both time series. For the reader's convenience, we denote the long-range dependence parameter $d$ by $H=d+\frac{1}{2}$ as it is common when dealing with fractional Gaussian noise and fractional Brownian motion. We refer to $H$ as \textit{Hurst parameter}, tracing back to the work of \cite{hurst:1951}.
\begin{figure}[h]
\centering
\includegraphics[width=1\textwidth]{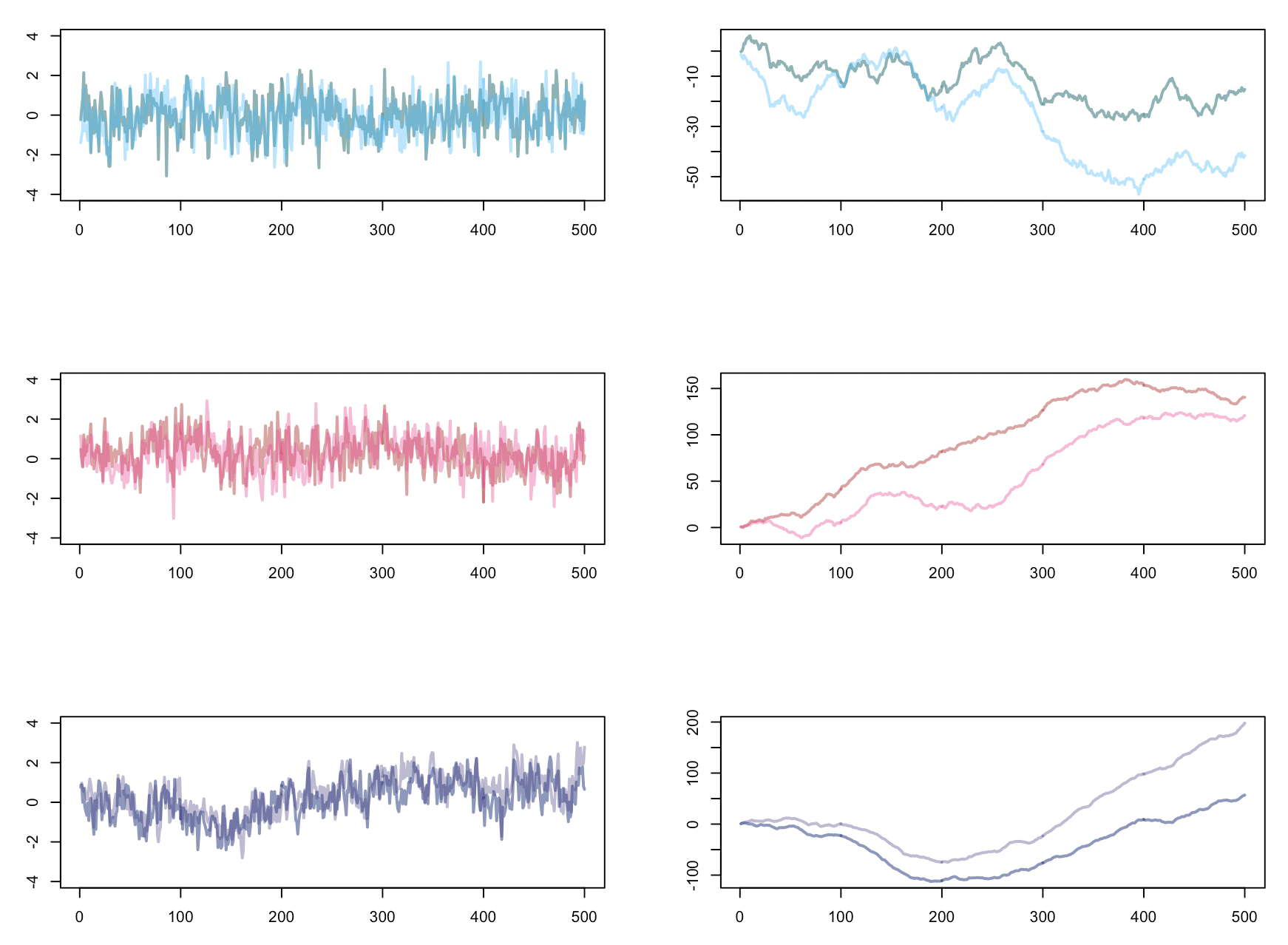}
\caption{Plots of $500$ data points of one path of two dependent fractional Gaussian noise processes (left) and the paths of the corresponding fractional Brownian motions (right) for different Hurst parameters: $H=0.7$ (top), $H=0.8$ (middle), $H=0.9$ (bottom).}
\label{fig: Timeseries}
\end{figure}~\newline
For $H=0.7$ and $H=0.8$ we generate $n=10^6$ samples, for $H=0.9$, we choose $n=2\cdot10^6$. We denote the correlation function of univariate fractional Gaussian noise by $r_H^{(1,1)}(k)$, $k\geq 0$. Then, we obtain $\left(Y_j^{(1)},Y_j^{(2)}\right)_j$ for $j=1,\ldots,n$:
\begin{align}
& Y_j^{(1)}=U_j^{(1)},\notag\\
& Y_j^{(2)}=\psi U_j^{(1)}+\phi U_j^{(2)},\label{eq:parametersimstudy}
\end{align}
for $\psi,\phi\in\mathbb{R}$.\newline
Note that this yields the following properties for the cross-correlations of the two processes for $k\geq 0$:
\begin{align*}
r_H^{(1,2)}(k)&=\mathbb{E}\left(Y_j^{(1)}Y_{j+k}^{(2)}\right)=\psi r_H^{(1,1)}(k)\\
r_H^{(2,1)}(k)&=r^{(1,2)}(-k)=\psi r_H^{(1,1)}(k)\\
r_H^{(2,2)}(k)&=\mathbb{E}\left(Y_j^{(2)}Y_{j+k}^{(2)}\right)=\left(\psi^2+\phi^2\right) r_H^{(1,1)}(k).
\end{align*}
We use $\psi=0.6$ and $\phi=0.8$ to get unit variance in the second process. Note that we choose the same Hurst parameter in both processes to get a better simulation result. The simulations of the processes $\left(Y_j^{(1)}\right)_{j\in\mathbb{Z}}$ and $\left(Y_j^{(2)}\right)_{j\in\mathbb{Z}}$ are visualized in Figure \ref{fig: Timeseries}. 
On the left-hand side the different fractional Gaussian noises depending on the Hurst parameter $H$ are displayed. They represent the stationary long-range dependent Gaussian \textit{increment} processes we need in the view of the limit theorems we derived in Section \ref{sec: OPD}. The processes in which we are comparing the coincident ordinal patterns, namely $\left(X_j^{(1)}\right)_{j\in\mathbb{Z}}$ and $\left(X_j^{(2)}\right)_{j\in\mathbb{Z}}$, are shown on the right-hand side in Figure \ref{fig: Timeseries}. The long-range dependent behaviour of the increment processes is very illustrative in these processes: roughly speaking they get smoother the larger the Hurst parameter gets.\newline
We turn to the simulation results for the asymptotic distribution of the estimator $\hat{p}_n$.
\begin{figure}[h]
\centering
\includegraphics[width=0.8\textwidth]{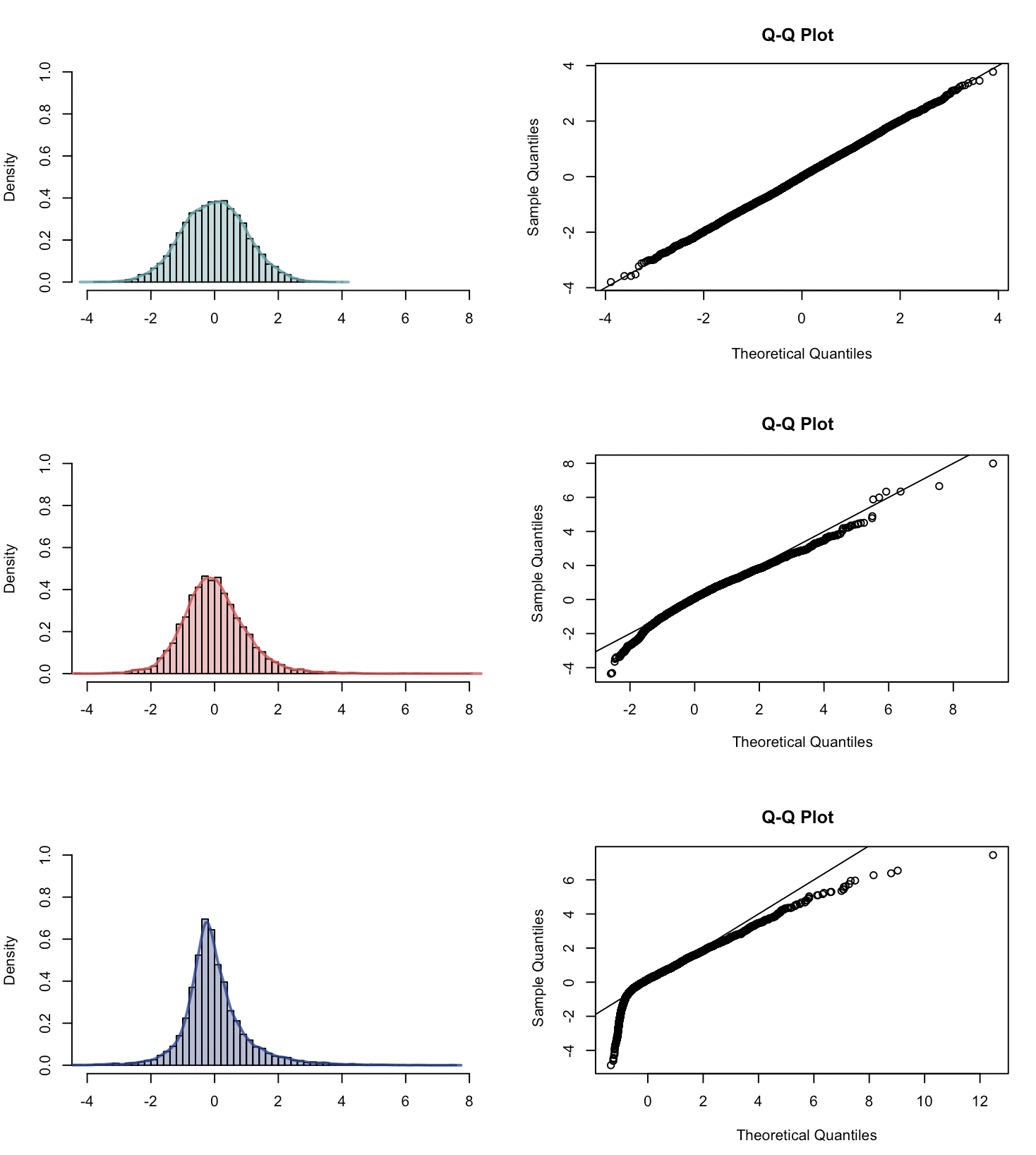}
\caption{Histogram, kernel density estimation and Q-Q plot with respect to the normal distribution ($H=0.7$) or to the Rosenblatt distribution of $\hat{p}_n-p$ with $h=2$ for different Hurst parameters: $H=0.7$ (top), $H=0.8$ (middle), $H=0.9$ (bottom).}
\label{fig: OPD}
\end{figure}
The first limit theorem is given in Theorem \ref{LT:OPDp} for $H=0.8$ and $H=0.9$. In the case $H=0.7$ a different limit theorem holds, see Theorem \ref{LT: SRDp}. Therefore, we turn to the simulation results of the asymptotic distribution of the estimator $\hat{p}_n$ of $p$, as shown in Figure \ref{fig: OPD} for pattern length $h=2$.
The asymptotic normality in case $H=0.7$ can be clearly observed. We turn to the interpretation of the simulation results of the distribution of $\hat{p}_n-p$ for $H=0.8$ and $H=0.9$ as weighted sum of the sample (cross-)correlations: we observe in the Q-Q plot for $H=0.8$ that the samples in the upper and lower tail deviate from the reference line. For $H=0.9$ a similar behaviour in the Q-Q plot is observed.\newline
We want to verify the result in Theorem \ref{LT:OPDpAlternative}, that it is possible by a different weighting, to express the limit distribution of $\hat{p}_n-p$ as the distribution of the sum of two independent standard Rosenblatt random variables. The simulated convergence result is provided in Figure \ref{fig: OPDalt}.
~\newline
We observe the standard Rosenblatt distribution.
\section{Conclusion and outlook}
We considered limit theorems in the context of estimation of ordinal pattern dependence in the long-range dependence setting. Pure long-range dependence as well as mixed cases of short- and long-range dependence were considered alongside with the transformed short-range dependent case. Therefore, we complemented the asymptotic results in \cite{schnurr:dehling:2017}. Hence, we made ordinal pattern dependence applicapable for long-range dependent data sets as they arise in the context of neurology, see \cite{karlekar:2014} or artifical intelligence, see \cite{shu:1999}. As these kinds of data were already investigated using ordinal patterns, see for example \cite{keller:sinn:lauffer2007}, this emphasizes the large practical impact of the ordinal approach in analyzing the dependence structure multivariate time series. This yields various research opportunities in these fields in the future.\newline\newline
This research was supported by the German Research Foundation (DFG) through the project \emph{ Ordinal-Pattern-Dependence: Grenzwerts\"atze und Strukturbr\"uche im langzeitabh\"angigen Fall mit Anwendungen in Hydrologie, Medizin und Finanzmathematik} (SCHN 1231/3-2).
\bibliographystyle{elsarticle-harv} 
\bibliography{OPDincontextLRD.bib}
\end{document}


\setcounter{section}{0}%
\appendix
    \renewcommand\thesection{A}%
\noindent All proofs in this appendix are taken from \cite{nuessgen:2021}, Chapter 3.
\section{Preliminary results}
\noindent Before turning to limit theorems, we introduce a possibility to decompose the $d$-dimensional Gaussian process $\left(Y_j\right)_{j\in\mathbb{Z}}$ using the Cholesky decomposition, see \cite{golub:2013}. Based on the definition of the multivariate normal distribution, see \cite{brockwell:davis:1991}, Definition 1.6.1, we find an upper triangular matrix $\tilde{A}$, such that $\tilde{A}\tilde{A}^t=\Sigma_d$. Then it holds, that
\begin{align}
Y_j\overset{\mathcal{D}}{=} \tilde{A} U^*_j,\label{eq:Choleskynormal1}
\end{align}
where $U^*_j$ is a $d$-dimensional Gaussian process where each $U_j^*$ has independent and identically $\mathcal{N}(0,1)$ distributed entries. We want to assure that $\left(U^*_j\right)_{j\in\mathbb{Z}}$ preserves the long-range dependent structure of $\left(Y_j\right)_{j\in\mathbb{Z}}$. Since we know from (2), that
\begin{align*}
\mathbb{E}\left(Y_j Y_{j+k}\right)=\Gamma_Y(k)\simeq k^{D-\frac{1}{2}I_d} L k^{D-\frac{1}{2}I_d}\quad (k\rightarrow\infty),
\end{align*}
the process $\left(U_j^*\right)$ has to fulfill
\begin{align}
\mathbb{E}\left(U^*_j U^*_{j+k}\right)=\Gamma_{U^*}(k)\simeq k^{D-\frac{1}{2}I_d}L_Uk^{D-\frac{1}{2}I_d}\quad (k\rightarrow\infty),\label{eq:LRDU}
\end{align}
with $L=\tilde{A}L_{U^*}\tilde{A}^t$.\newline
Then it holds for all $n\in\mathbb{N}$ that
\begin{align}
\left(Y_j,~j=1,\ldots,n\right)\overset{\mathcal{D}}{=} \left(\tilde{A}U_j^*,~j=1,\ldots,n\right).\label{eq:Choleskynormal2}
\end{align}\quad\newline
Note that the assumption in \eqref{eq:LRDU} is only well-defined because we assumed $\left|r^{(p,q)}(k)\right|<1$ for $k\geq 1$ and $p,q=1,\ldots,d$ in (1). This becomes clear in the following considerations. In the proofs of the theorems in this chapter, we do not only need a decomposition of $Y_j$, but also of $Y_{j,h}$. As $Y_{j,h}$ is still a multivariate Gaussian process, the covariance matrix of $Y_{j,h}$ given by $\Sigma_{d,h}$ is positive definite. Hence, it is possible to find a upper triangular matrix $A$, such that $AA^t=\Sigma_{d,h}$. It holds that
\begin{align*}
Y_{j,h}\overset{\mathcal{D}}{=}AU_{j,h}
\end{align*}
for 
\begin{align*}
U_{j,h}=\left(U_{(j-1)h+1}^{(1)},\ldots,U_{jh}^{(1)},\ldots,U_{(j-1)h+1}^{(d)},\ldots,U_{jh}^{(d)}\right)^t.
\end{align*}
The random vector $U_{j,h}$ consists of $(d\cdot h)$ independent and standard normally distributed random variables. We notice the different structure of $U_{j,h}$ compared to $Y_{j,h}$. We assure that for consecutive $j$ the entries in $U_{j,h}$ are all different while there are identical entries, for example in $Y_{1,h}=\left(Y_1^{(1)},Y_2^{(1)},\ldots,Y_h^{(d)}\right)^t$ and $Y_{2,h}=\left(Y_2^{(1)},\ldots,Y_h^{(d)},Y_{h+1}^{(d)}\right)^t$. This complicates our aim that 
\begin{align}
\left(Y_{j,h},j=1,\ldots,n\right)^t\overset{\mathcal{D}}{=}\left(AU_{j,h},~j=1,\ldots,n\right)^t\label{eq:jointconvergence}
\end{align}
holds.\newline 
The special structure of $\left(Y_{j,h}\right)_{j\in\mathbb{Z}}$, namely, that it is consisting of $h$ consecutive entries of each marginal process $\left(Y_j^{(p)}\right)$, $p=1,\ldots,d$, alongside with the dependence between two random vectors in the process $\left(Y_{j,h}\right)$, has to be reflected in the covariance matrix of $\left(U_{j,h},~j=1,\ldots,n\right)$. Hence, we need to check whether such a vector $\left(U_{j,h},~j=1,\ldots,n\right)$ exists, i.e., if there is a positiv semi-definite matrix that fulfills these conditions. We define $\textbf{A}$ as a block diagonal matrix with $A$ as main-diagonal blocks and all off-diagonal blocks as $dh\times dh$-zero matrix.
We denote the covariance matrix of $\left(Y_{j,h},j=1,\ldots,n\right)^t$ by $\Sigma_{Y,n}$ and define the following matrix:
\begin{align}
\Sigma_{U,n}:=\mathrm{inv}\left(\textbf{A}\right)\Sigma_{Y,n}\mathrm{inv}\left(\textbf{A}^t\right).\label{eq:solution}
\end{align}
We know that $\Sigma_{Y,n}$ is positive semi-definite for all $n\in\mathbb{N}$ because $\left(Y_j\right)$ is a Gaussian process. Mathematically described that means
\begin{align}
x^t \Sigma_{Y,n} x \geq 0,\label{eq:possemidef}
\end{align}
for all $x=\left(x_1,\ldots,x_{nhd}\right)^t\in\mathbb{R}^{nhd}$. We conclude
\begin{align*}
x^t \Sigma_{U,n} x &\:\:\,= x^t\mathrm{inv}\left(\textbf{A}\right)\Sigma_{Y,n}\mathrm{inv}\left(\textbf{A}^t\right) x\\
&\:\:\,=\left(\mathrm{inv}\left(\textbf{A}^t\right)x\right)^t \Sigma_{Y,n} \left(x^t\mathrm{inv}\left(\textbf{A}\right)\right)\\
&\overset{\eqref{eq:possemidef}}{\geq} 0.
\end{align*}
Therefore, $\Sigma_{U,n}$ is a positive semi-definite matrix for all $n\in\mathbb{N}$ and the random vector 
\begin{align*}
\left(U_{j,h},~j=1,\ldots,n\right)^t\mathcal{N}\sim\left(0,\Sigma_{U,n}\right)
\end{align*}
exists and \eqref{eq:jointconvergence} holds. Note that we do not have any further information on the dependence structure within the process $\left(U_j\right)$, in general, this process does neither exhibit long-range dependence nor independence nor stationarity.

We continue with two preparatory results that are as well necessary for proving Theorem 2.1.\newline
\begin{Lem}\label{lemm:HR2preceding}
Let $\left(Y_j\right)_{j\in\mathbb{Z}}$ be a $d$-dimensional Gaussian process as defined in (1) that fulfills (2) with $d_1=\ldots=d_d=d^*$, such that
\begin{align*}
\Gamma_Y(k)=\mathbb{E}\left(Y_jY_{j+k}^t\right)\simeq L k^{2d^*-1},\quad(k\rightarrow\infty).
\end{align*}
Let $C_2$ be a normalization constant,
\begin{align*}
C_2=\frac{1}{2d^*\left(4d^*-1\right)}
\end{align*}
and let $B_Y$ be an upper triangular matrix, such that 
\begin{align*}
B_YB_Y^t=L.
\end{align*}
Further, for $l\in\mathbb{N}$ we have
\begin{align*}
\hat{\Gamma}_{Y,n}(l)=\frac{1}{n-l}\sum_{j=1}^{n-l} Y_j Y_{j+l}^t.
\end{align*}
Then, for $h\in\mathbb{N}$ it holds that
\begin{align*}
&\left(n^{1-2d^*}\left(C_2\right)^{-1/2} \left(B_Y\otimes B_Y \right)^{-1}   \mathrm{vec}\left(\hat{\Gamma}_{n}(l)-\Gamma(l)\right),~l=0,\ldots,h-1\right)\\
& \quad\quad \quad\quad \quad\quad \quad\quad \quad\quad \quad\quad\xrightarrow{D} \left(\mathrm{vec}\left(Z^{(p,q)}_{2,d^*+1/2}(1)\right)_{p,q=1,\ldots,d},~l=0,\ldots,h-1\right),
\end{align*}
where $Z^{(p,q)}_{2,d^*+1/2}(1)$ has the spectral domain representation
\begin{align*}
Z^{(p,q)}_{2,d^*+1/2}(1)=K_{p,q}(d^*)\int_{\mathbb{R}^2}^{\prime\prime} \frac{\exp\left(i\left(\lambda_1+\lambda_2\right)\right)-1}{i\left(\lambda_1+\lambda_2\right)}\left|\lambda_1\lambda_2  \right|^{-d^*}\tilde{B}^{(p)}\left( \mathrm{d}\lambda_1\right)\tilde{B}^{(q)}\left( \mathrm{d}\lambda_2\right)
\end{align*}
where 
\begin{align*}
K_{p,q}^2(d^*)=\begin{cases} &\frac{1}{2C_2\left( 2\Gamma(1-2d^*)\sin\left(\pi d^*\right) \right)^2},\quad p=q\\
&\frac{1}{C_2\left( 2\Gamma(1-2d^*)\sin\left(\pi d^*\right) \right)^2},\:\,\quad p\neq q.
\end{cases}
\end{align*} and
 $\tilde{B}(\mathrm{d}\lambda)=\left(\tilde{B}^{(1)}(\mathrm{d}\lambda),\ldots,\tilde{B}^{(d)}(\mathrm{d}\lambda)\right)$ is a multivariate Hermitian-Gaussian random measure as defined in (12).
\end{Lem}
\begin{pf}
First, we can use \eqref{eq:Choleskynormal1}:
\begin{align*}
Y_j \overset{\mathcal{D}}{=} \tilde{A} U_j^*,
\end{align*}
such that $\left(U_j^*\right)$ is a multivariate Gaussian process with $U^*_j\sim\mathcal{N}\left(0,I_d\right)$ and $\left(U_j^*\right)$ is still long-range dependent, see \eqref{eq:LRDU}.
It is possible to decompose the sample cross-covariance matrix $\hat{\Gamma}_{Y,n}(l)-\Gamma_Y(l)$ with respect to $\left(Y_j\right)$ at lag $l$ given by
\begin{align*}
\hat{\Gamma}_{Y,n}(l)-\Gamma_Y(l)=\frac{1}{n-l}\sum_{j=1}^{n-l} Y_j Y_{j+l}^t - \mathbb{E}\left(Y_jY_{j+l}^t\right)
\end{align*}
to 
\begin{align*}
\hat{\Gamma}_{Y,n}(l)-\Gamma_Y(l)\overset{\mathcal{D}}{=}   \tilde{A}\left( \hat{\Gamma}_{U^*,n}(l)-\Gamma_{U^*}(l)\right) \tilde{A}^t,
\end{align*}
where we define the sample cross-covariance matrix $\hat{\Gamma}_{U^*,n}(l)-\Gamma_{U^*}(l)$ with respect to $\left(U_j^*\right)$ at lag $l$ by
\begin{align*}
 \hat{\Gamma}_{U^*,n}(l)-\Gamma_{U^*}(l)=\frac{1}{n-l}\sum_{j=1}^{n-l} U^*_j U^*_{j+l} - \mathbb{E}\left(U_j^* U_{j+l}^* \right).
\end{align*}
Each entry of 
\begin{align*}
\hat{\Gamma}_{U^*,n}(l)-\Gamma_{U^*}(l)=\left(\hat{r}_{n,U^*}^{(p,q)}(l)-r_{U^*}^{(p,q)}(l)\right)_{p,q=1,\ldots,d}
\end{align*}
is given by
\begin{align*}
\hat{r}_{n,U^*}^{(p,q)}(l)-r_{U^*}^{(p,q)}(l):=\sum_{j=1}^n U_j^{*~(p)}U_{j+l}^{*~(q)}-\mathbb{E}\left(U_j^{*~(p)}U_{j+l}^{*~(q)}\right).
\end{align*}
Following \cite{dueker:2017}, proof of Lemma 7.4, the limit distribution of
\begin{align*}
\left( \hat{\Gamma}_{U^*,n}(l)-\Gamma_{U^*}(l),~l=0,\ldots,h-1\right)
\end{align*}
is equal to the limit distribution of
\begin{align*}
\left( \hat{\Gamma}_{U^*,n}(0)-\Gamma_{U^*}(0),~l=0,\ldots,h-1\right).
\end{align*}
We recall the assumption that $d^*=d_p$ for all $p=1,\ldots,d$.
We follow \cite{arcones:1994}, Theorem 6 and use the Cramer-Wold device: Let $a_{1,1},a_{1,2},\ldots, a_{d,d}\in\mathbb{R}$. We are interested in the asymptotic behaviour of
\begin{align*}
&n^{1-2d^*}\sum_{p,q=1}^d a_{p,q}\left(\hat{r}_{n,U}^{(p,q)}(0)-r^{(p,q)}_U(0)\right)\\
=&n^{-2d^*} \sum_{j=1}^n \sum_{p,q=1}^d a_{p,q} \left(  U_j^{*~(p)}U_{j}^{*~(q)}-\mathbb{E}\left(U_j^{*~(p)}U_{j}^{*~(q)}\right)\right).
\end{align*}
We consider the function 
\begin{align}
f\left(U_j^*\right)=\sum_{p,q=1}^d a_{p,q} \left(  U_j^{*~(p)}U_{j}^{*~(q)}-\mathbb{E}\left(U_j^{*~(p)}U_{j}^{*~(q)}\right) \right)\label{eq:baldistsfertig}
\end{align}
and may apply Theorem 6 in \cite{arcones:1994}. Using the Hermite decomposition of $f$ as given in (11), we observe that $f$ and, therefore, $a_{p,q}$, $p,q=1,\ldots,d$, only affects the Hermite coefficients. Indeed, using \cite{beran:feng:ghosh:kulik:2013}, Lemma 3.5, the Hermite coefficients reduce to $a_{p,q}$ for each summand on the right-hand side in \eqref{eq:baldistsfertig}. Hence, we can state
\begin{align}
n^{-2d^*} \sum_{j=1}^n \sum_{p,q=1}^d a_{p,q} &\left(  U_j^{*~(p)}U_{j}^{*~(q)}-\mathbb{E}\left(U_j^{*~(p)}U_{j}^{*~(q)}\right)\right)\\
&\xrightarrow{\mathcal{D}} \sum_{p,q=1}^d a_{p,q} Z^{(p,q)}_{2,d^*+1/2}(1),\label{eq:petermajor}
\end{align}
where $Z^{(p,q)}_{2,d^*+1/2}(1)$ has the spectral domain representation
\begin{align}
Z^{(p,q)}_{2,d^*+1/2}(1)=K_{p,q}(d^*)\int_{\mathbb{R}^2}^{\prime\prime} \frac{\exp\left(i\left(\lambda_1+\lambda_2\right)\right)-1}{i\left(\lambda_1+\lambda_2\right)}\left|\lambda_1\lambda_2  \right|^{-d^*} \tilde{B}^{(p)}\left(\mathrm{d}\lambda_1\right) \tilde{B}^{(q)}\left(\mathrm{d}\lambda_2\right)\label{eq:hermiteprocess}
\end{align}
where
\begin{align*}
K^2_{p,q}(d^*)=\begin{cases} &\frac{1}{2C_2\left( 2\Gamma(1-2d^*)\sin\left(\pi d^*\right) \right)^2},\quad p=q\\
&\frac{1}{C_2\left( 2\Gamma(1-2d^*)\sin\left(\pi d^*\right) \right)^2},\quad\quad p\neq q.
\end{cases}
\end{align*}
and $\tilde{B}(\mathrm{d}\lambda)=\left(\tilde{B}^{(1)}(\mathrm{d}\lambda),\ldots,\tilde{B}^{(d)}(\mathrm{d}\lambda)\right)$ is an appropriate multivariate Hermitian-Gaussian random measure.
Thus, we proved convergence in distribution of the sample-cross correlation matrix:
\begin{align*}
n^{1-2d^*}\left(\hat{\Gamma}_{U^*,n}(0)-\Gamma_{U^*}(0)\right)\xrightarrow{D} \left(Z^{(p,q)}_{2,d^*+1/2}(1)\right)_{p,q=1,\ldots,d}.
\end{align*}
We take a closer look at the covariance matrix of $\mathrm{vec}\left(\hat{\Gamma}_{U^*,n}(0)-\Gamma_{U^*}(0)\right)$. Following \cite{dueker:2019}, Lemma 5.7, we observe
\begin{align*}
&n^{1-2d^*}\left(4d^*\left(4d^*-1\right)\right)^{1/2}\mathrm{Cov}\left(\mathrm{vec}\left(\hat{\Gamma}_{U^*,n}(0)-\Gamma_{U^*}(0)\right),\mathrm{vec}\left(\hat{\Gamma}_{U^*,n}(0)-\Gamma_{U^*}(0)\right)\right)\\
=&\left(I_{d^2}+K_{d^2}\right)\left(L_{U^*} \otimes L_{U^*}\right),
\end{align*}
with $L_{U^*}$ as defined in \eqref{eq:Choleskynormal2} and $\otimes$ denotes the Kronecker product. Furthermore $K_{d}$ denotes the commutation matrix that that transforms $\mathrm{vec}(A)$ into $\mathrm{vec}\left(A^t\right)$ for $A\in\mathbb{R}^{d\times d}$. For details see \cite{magnus:1979}.\newline
Hence, the covariance matrix of the vector of the sample cross-covariances is fully specified by the knowledge of $L_{U^*}$ as it arises in the context of long-range dependence in \eqref{eq:Choleskynormal2}.\newline
We obtain a relation between $L$ and $L_{U^*}$, since
\begin{align*}
\Gamma_Y(\cdot)=\tilde{A}\Gamma_U(\cdot)\tilde{A}^t.
\end{align*}
Both
\begin{align*}
\Gamma_Y(k)\simeq  Lk^{2d^*-1}~(k\rightarrow\infty)
\end{align*}
and
\begin{align*}
\Gamma_{U^*}(k)\simeq L_{U^*}k^{2d^*-1}~(k\rightarrow\infty)
\end{align*}
hold and we obtain
\begin{align*}
L=\tilde{A}L_{U^*}\tilde{A}^t.
\end{align*}
We study the covariance matrix of $\mathrm{vec}\left(\hat{\Gamma}_{Y,n}(0)-\Gamma_{Y}(0)\right)$:
\begin{align}
&n^{1-2d^*}\left(4d^*\left(4d^*-1\right)\right)^{1/2}\mathrm{Cov}\left(\mathrm{vec}\left(\hat{\Gamma}_{Y,n}(0)-\Gamma_{Y}(0)\right),\mathrm{vec}\left(\hat{\Gamma}_{Y,n}(0)-\Gamma_{Y}(0)\right)^t\right)\notag\\
\rightarrow &\quad \left(I_{d^2}+K_{d^2}\right)\left( L \otimes L\right) \label{eq:covarianceRosenblatt} \\
=&\quad\left(I_{d^2}+K_{d^2}\right) \left(\tilde{A}L_{U^*}\tilde{A}^t\right)\otimes \left(\tilde{A}L_{U^*}\tilde{A}^t\right)\notag\\
=&\quad \left(I_{d^2}+K_{d^2}\right)\left(\tilde{A}\otimes\tilde{A}\right)\cdot\left(L_{U^*}\otimes L_{U^*} \right)\cdot \left(\tilde{A}^t\otimes\tilde{A}^t\right).\notag
\end{align}
Let $B_{U^*}$ be an upper triangular matrix, such that 
\begin{align*}
B_{U^*}B_{U^*}^t:=L_{U^*}.
\end{align*}
We know that such a matrix exists because $L_{U^*}$ is positive definite.
Analogously, we define $B_{Y}$:
\begin{align*}
B_Y:=\tilde{A} B_{U^*}.
\end{align*}
Then, it holds that
\begin{align*}
B_YB_Y^t=L.
\end{align*}
We arrive at
\begin{align*}
& n^{1-2d^*}\left(C_2\right)^{-1/2} \left(B_{Y}\otimes B_{Y}\right)^{-1}  \mathrm{vec}\left(\hat{\Gamma}_{Y,n}(0)-\Gamma_Y(0)\right)\\
&\overset{\mathcal{D}}{=}n^{1-2d^*}\left(C_2\right)^{-1/2} \left(B_{U^*}\otimes B_{U^*}\right)^{-1}\left(A\otimes A\right)^{-1} \mathrm{vec}\left( \tilde{A} \left(\hat{\Gamma}_{U^*,n}(0)-\Gamma_{U^*}(0)\right) \tilde{A}^t\right) \\
&=n^{1-2d^*}\left(C_2\right)^{-1/2} \left(B_{U^*}\otimes B_{U^*}\right)^{-1}\mathrm{vec}\left(\hat{\Gamma}_{U^*,n}(0)-\Gamma_{U^*}(0)\right)\\
& \xrightarrow{D} \mathrm{vec}\left(Z^{(p,q)}_{2,d^*+1/2}(1)\right)_{p,q=1,\ldots,d},
\end{align*}
where $Z^{(p,q)}_{2,d^*+1/2}(1)$ has the spectral domain representation
\begin{align*}
Z^{(p,q)}_{2,d^*+1/2}(1)=K_{p,q}(d^*)\int_{\mathbb{R}^2}^{\prime\prime} \frac{\exp\left(i\left(\lambda_1+\lambda_2\right)\right)-1}{i\left(\lambda_1+\lambda_2\right)}\left|\lambda_1\lambda_2  \right|^{-d^*}\tilde{B}^{(p)}\left( \mathrm{d}\lambda_1\right)\tilde{B}^{(q)}\left( \mathrm{d}\lambda_2\right)
\end{align*}
where 
\begin{align*}
K^2_{p,q}(d^*)=\begin{cases} &\frac{1}{2C_2\left( 2\Gamma(1-2d^*)\sin\left(\pi d^*\right) \right)^2},\quad p=q\\
&\frac{1}{C_2\left( 2\Gamma(1-2d^*)\sin\left(\pi d^*\right) \right)^2},\:\,\quad p\neq q.
\end{cases}
\end{align*}
and $\tilde{B}(\mathrm{d}\lambda)=\left(\tilde{B}^{(1)}(\mathrm{d}\lambda),\ldots,\tilde{B}^{(d)}(\mathrm{d}\lambda)\right)$ is a multivariate Hermitian-Gaussian random measure as defined in (12). Note that the standardization on the left-hand side is appropriate since the covariance matrix of $\mathrm{vec}\left(Z_{2,d^*+1/2}(1)\right)$ is given by
\begin{align}
\mathbb{E}\Bigg (K^2(d^*)\int_{\mathbb{R}^2}^{\prime\prime}&\int_{\mathbb{R}^2}^{\prime\prime}E_{\lambda_1,\lambda_2}\overline{E_{\lambda_3,\lambda_4}}   \mathrm{vec}\left(\tilde{B}\left(\mathrm{d}\lambda_1\right)\left(\tilde{B}\left( \mathrm{d}\lambda_2\right)\right)^t \right)  \notag\\
&\quad\quad\quad\quad\quad\quad\quad\quad\quad \left(   \mathrm{vec} \left(\overline{\tilde{B}\left(\mathrm{d}\lambda_3\right)\left( {\tilde{B}}\left(\mathrm{d}\lambda_4\right)\right)^t }\right)    \right)^t      \Bigg).\label{eq:excludediag}
\end{align}
by denoting
\begin{align*}
E_{\lambda_1,\lambda_2}:= \frac{\exp\left(i\left(\lambda_1+\lambda_2\right)\right)-1}{i\left(\lambda_1+\lambda_2\right)}\left|\lambda_1\lambda_2  \right|^{-d^*}.
\end{align*}
We observe
\begin{align}
&\mathbb{E}\left(\mathrm{vec}\left(\tilde{B}\left(\mathrm{d}\lambda_1\right)\tilde{B}\left(\mathrm{d}\lambda_2\right)^t  \right) \left( \mathrm{vec} \left(\overline{\tilde{B}\left(\mathrm{d}\lambda_3\right)\left( {\tilde{B}}\left(\mathrm{d}\lambda_4\right)\right)^t }\right)     \right)^t   \right)\notag\\
=&\begin{cases}
&I_{d^2}\mathrm{d}\lambda_1\mathrm{d}\lambda_2,\quad \left|\lambda_1\right|=\left|\lambda_3\right|\neq \left|\lambda_2\right|=\left|\lambda_4\right|,\\
& K_{d^2}\mathrm{d}\lambda_1\mathrm{d}\lambda_2,\quad \left|\lambda_1\right|=\left|\lambda_4\right|\neq \left|\lambda_2\right|=\left|\lambda_3\right|,
\end{cases}\label{eq:considerationssoimportant}
\end{align}
following \cite{dueker:2019}, (27). Neither the case $\left|\lambda_1\right|=\left|\lambda_2\right|$ nor $ \left|\lambda_3\right|=\left|\lambda_4\right|$ has to be incorporated as the diagonals are excluded in the integration in \eqref{eq:excludediag}.
 $\hfill\Box$\newline
\end{pf}
\begin{Kor}\label{korl:precedingHR2}
Under the assumptions of Lemma \ref{lemm:HR2preceding}, there is a different representation of the limit random vector.
For $h\in\mathbb{N}$ we obtain
\begin{align*}
\left(n^{1-2d^*}\left(C_2\right)^{-1/2} \mathrm{vec}\left(\hat{\Gamma}_{n}(l)-\Gamma(l)\right),~l=0,\ldots,h-1\right)\xrightarrow{D} \left(\mathrm{vec}\left(Z_{2,d^*+1/2}(1)\right)~l=0,\ldots,h-1\right),
\end{align*}
where $\mathrm{vec}\left(Z_{2,d^*+1/2}(1)\right)$ has the spectral domain representation
\begin{align*}
\mathrm{vec}\left(Z_{2,d^*+1/2}(1)\right)=D_{K\left(d^*\right)}\int_{\mathbb{R}^2}^{\prime\prime} \frac{\exp\left(i\left(\lambda_1+\lambda_2\right)\right)-1}{i\left(\lambda_1+\lambda_2\right)}\left|\lambda_1\lambda_2  \right|^{-d^*}\mathrm{vec}\left(\tilde{B}_L\left(\mathrm{d}\lambda_1\right)\tilde{B}_L\left( \mathrm{d}\lambda_2\right)^t\right).
\end{align*}
The matrix $D_{K\left(d^*\right)}$ is a diagonal matrix, 
\begin{align*}
D_{K\left(d^*\right)}=\mathrm{diag}\left(\mathrm{vec}\left(K\left(d^*\right)\right)\right),
\end{align*}
and $K(d^*)=\left(K_{p,q}(d^*)\right)_{p,q=1,\ldots,d}$ is such that
\begin{align*}
K^2(d^*)_{p,q}=\begin{cases} &\frac{1}{2C_2\left( 2\Gamma(1-2d^*)\sin\left(\pi d^*\right) \right)^2},\quad p=q\\
&\frac{1}{C_2\left( 2\Gamma(1-2d^*)\sin\left(\pi d^*\right) \right)^2},\quad\:\, p\neq q.
\end{cases}
\end{align*}
Furthermore, $\tilde{B}_L(\mathrm{d}\lambda)$ is a multivariate Hermitian-Gaussian random measure that fulfills 
\begin{align*}
\mathbb{E}\left(\tilde{B}_L(\mathrm{d}\lambda) \tilde{B}_L(\mathrm{d}\lambda)^*\right)=L~\mathrm{d}\lambda.
\end{align*}
\end{Kor}
\begin{pf}
The proof is an immediate consequence of Lemma \ref{lemm:HR2preceding} using $\tilde{B}_L(\mathrm{d}\lambda)=B_Y\tilde{B}(\mathrm{d}\lambda)$ with $B_YB_Y^t=L$ and $\tilde{B}(\mathrm{d}\lambda)$ as defined in (12). $\hfill\Box$\newline
\end{pf}
    \renewcommand\thesection{B}%
\section{Proof of Theorem 2.1}
\begin{pf}
Without loss of generality, we assume $\mathbb{E}\left( f\left( Y_{j,h}\right)  \right)=0$.
Following the argumentation in \cite{betken:2019}, Theorem 5.9, we first remark that $Y_{j,h}\overset{\mathcal{D}}{=}AU_{j,h}$ with $U_{j,h}$ and $A$ as described in \eqref{eq:jointconvergence} and \eqref{eq:solution}.
We now want to study the asymptotic behavior of the partial sum $\sum\limits_{j=1}^n f^*\left(U_{j}\right)$ where $f^*\left(U_{j,h}\right):=f\left(AU_{j,h}\right)\overset{\mathcal{D}}{=}f\left(Y_{j,h}\right)$. Since $m\left(f^*,I_{dh}\right)=m\left(f\circ A, I_{dh}\right)=m\left(f,\Sigma_{d,h}\right)=2$, see \cite{beran:feng:ghosh:kulik:2013}, Lemma 3.7, hence, we know by \cite{arcones:1994}, Theorem 6, that these partial sums are dominated by the second order terms in the Hermite expansion of $f^*$:
\begin{align*}
\sum_{j=1}^n f^*\left(U_{j,h}\right)&\sum_{j=1}^n \sum_{l_1+\ldots+l_{dh}=2}\mathbb{E}\left(f^*\left(U_{j,h}\right) H_{l_1,\ldots,l_{dh}}\left(U_{j,h}\right)\right)H_{l_1,\ldots,l_{dh}}\left(U_{j,h}\right)+o_{\mathbb{P}}\left(n^{2d^*}\right).
\end{align*}
This follows from the multivariate extension of the Reduction Theorem as proved in \cite{arcones:1994}. We obtain
\begin{align*}
&\sum_{l_1+\ldots+l_{dh}=2}\mathbb{E}\left(f^*\left(U_{j,h}\right) H_{l_1,\ldots,l_{dh}}\left(U_{j,h}\right)\right)H_{l_1,\ldots,l_{dh}}\left(U_{j,h}\right)\\
=&\sum_{i=1}^{dh} \mathbb{E}\left(f^*\left(U_{j,h}\right) \left(\left(U_{j,h}^{(i)}\right)^2-1\right) \right) \left(\left(U_{j,h}^{(i)}\right)^2-1\right)+\sum_{1\leq i,k\leq dh,i\neq k} \mathbb{E}\left(f^*\left(U_{j,h}\right)U_{j,h}^{(i)}U_{j,h}^{(k)}\right)U_{j,h}^{(i)}U_{j,h}^{(k)}\\
=&\sum_{i=1}^{dh} \mathbb{E}\left(f^*\left(U_{j,h}\right) \left(U_{j,h}^{(i)}\right)^2\right) \left(\left(U_{j,h}^{(i)}\right)^2-1\right)+\sum_{1\leq i,k\leq dh,i\neq k} \mathbb{E}\left(f^*\left(U_{j,h}\right)U_{j,h}^{(i)}U_{j,h}^{(k)}\right)U_{j,h}^{(i)}U_{j,h}^{(k)},
\end{align*}
since $\mathbb{E}\left(f^*\left(U_{j,h}\right)\right)=\mathbb{E}\left(f\left(Y_{j,h}\right)\right)=0$. This results in:
\begin{align}
&\sum_{i=1}^{dh} \mathbb{E}\left(f^*\left(U_{j,h}\right) \left(U_{j,h}^{(i)}\right)^2\right) \left(\left(U_{j,h}^{(i)}\right)^2-1\right)+\sum_{1\leq i,k\leq dh,i\neq k} \mathbb{E}\left(f^*\left(U_{j,h}\right)U_{j,h}^{(i)}U_{j,h}^{(k)}\right)U_{j,h}^{(i)}U_{j,h}^{(k)}\notag\\
=& \sum_{1\leq i,k\leq dh} \mathbb{E}\left( f^*\left(U_{j,h}\right) U_{j,h}^{(i)}  U_{j,h}^{(k)} \right)  U_{j,h}^{(i)}  U_{j,h}^{(k)} -\sum_{i=1}^{dh} \mathbb{E}\left(f^*\left(U_{j,h}\right) \left(U_{j,h}^{(i)}\right)^2\right).\label{eq: th2 umformung1}
\end{align}
Note that 
\begin{align}
 \sum_{1\leq i,k\leq dh} \mathbb{E}\left( f^*\left(U_{j,h}\right) U_{j,h}^{(i)}  U_{j,h}^{(k)} \right) \mathbb{E}\left( U_{j,h}^{(i)}  U_{j,h}^{(k)} \right) =\sum_{i=1}^{dh} \mathbb{E}\left(f^*\left(U_{j,h}\right) \left(U_{j,h}^{(i)}\right)^2\right) \label{th1:erwartungswertgleichheit}
\end{align}
since the entries of $U_{j,h}$ are independent for fixed $j$ and identically $\mathcal{N}(0,1)$ distributed. So the subtrahend in \eqref{eq: th2 umformung1} equals the expected value of the minuend.\newline
Define $B:=\left(b_{i,k}\right)_{1\leq i,k\leq dh}\in\mathbb{R}^{(dh)\times(dh)}$ with $b_{i,k}:=\mathbb{E}\left(f^*\left(U_{j,h}\right) U_{j,h}^{(i)} U_{j,h}^{(k)}\right)=\mathbb{E}\left(f^*\left(U_{1}\right) U_{1}^{(i)} U_{1}^{(k)}\right)$ since we are considering a stationary process. We obtain
\begin{align*}
B=\mathbb{E}\left(U_{j,h} f^*\left(U_{j,h}\right) U_{j,h}^t\right)=\mathbb{E}\left(A^{-1}Y_{j,h}f\left(Y_{j,h}\right) Y_{j,h}^t \left(A^{-1}\right)^t \right).
\end{align*} \newline
Hence, we can state the following:
\begin{align}
\sum_{1\leq i,k\leq dh} \mathbb{E}\left(f^*\left(U_{j,h}\right)U_{j,h}^{(i)}U_{j,h}^{(k)}\right) U_{j,h}^{(i)}U_{j,h}^{(k)} &= U_{j,h}^t B U_{j,h}\notag\\
& \overset{\mathcal{D}}{=}Y_{j,h}^t\left( A^{-1}\right)^t B A^{-1}Y_{j,h}\notag \\
&= Y_{j,h}^t\left( A^{-1}\right)^t A^{-1}\mathbb{E}\left(Y_{j,h}f\left(Y_{j,h}\right) Y_{j,h}^t\right) \left(A^{-1}\right)^t A^{-1}Y_{j,h} \notag\\
&=Y_{j,h}^t \Sigma_{d,h}^{-1} \mathbb{E}\left(Y_{j,h}f\left(Y_{j,h}\right) Y_{j,h}^t\right) \Sigma_{d,h}^{-1} Y_{j,h}\notag\\
&=Y_{j,h}^t \mathbb{A} Y_{j,h}\notag\\
&=\sum_{1\leq i,k\leq dh} Y_j^{(i)}Y_j^{(k)}\alpha_{ik}, \label{th1: umformung1}
\end{align}
where we define $\mathbb{A}:=\left(\alpha_{ik}\right)_{1\leq ik\leq dh}:=\Sigma_{d,h}^{-1} C \Sigma_{d,h}^{-1}$, with $C:= \mathbb{E}\left(Y_{j,h}f\left(Y_{j,h}\right) Y_{j,h}^t\right) $ as the matrix of second order Hermite coefficients, in contrast to $B$ now with respect to the original considered process $\left(Y_{j,h}\right)_{j\in\mathbb{Z}}$. \newline\newline
Remembering the special structure of $Y_{j,h}=\left(Y_j^{(1)},\ldots,Y_{j+h-1}^{(1)},\ldots,Y_j^{(d)},Y_{j+h-1}^{(d)}\right)^t$, namely that $Y_{j,h}^{(k)}=Y_{j+(k\mod h)-1}^{\left(\lfloor \frac{k-1}{h} \rfloor +1\right)}$, $k=1,\ldots,dh$ we can see that
\begin{align}
\sum_{j=1}^n \sum_{1\leq i k \leq dh} Y_{j,h}^{(i)}Y_{j,h}^{(k)}\alpha_{ik}
&=\sum_{j=1}^n \sum_{1\leq i k \leq dh} Y_{j+(i\mod h)-1}^{\left(\left\lfloor \frac{i-1}{h} \right\rfloor +1\right)}Y_{j+(k\mod h)-1}^{\left(\left\lfloor \frac{k-1}{h} \right\rfloor +1\right)}\alpha_{ik} \notag\\
&=\sum_{j=1}^n \sum_{p,q=1}^d \sum_{i,k=1}^h Y_{j+i-1}^{(p)} Y_{j+k-1}^{(q)} \alpha_{ik}^{(p,q)},\label{eq: break down original process}
\end{align}
where we divide
\begin{align*}
\mathbb{A}=\begin{pmatrix}
\mathbb{A}^{(1,1)} && \mathbb{A}^{(1,2)} && \ldots && \mathbb{A}^{(1,d)}\\
\mathbb{A}^{(2,1)} && \mathbb{A}^{(2,2)} && \ldots && \mathbb{A}^{(2,d)}\\
\vdots && \vdots && \vdots \\
\mathbb{A}^{(d,1)} && \mathbb{A}^{(d,2)} && \ldots && \mathbb{A}^{(d,d)}\\
\end{pmatrix},
\end{align*}
with $\mathbb{A}^{(p,q)}=\left(\alpha^{(p,q)}_{i,k}\right)_{1\leq i,k\leq h}\in\mathbb{R}^{h\times h}$ such that $\alpha^{(p,q)}_{i,k}=\alpha_{i+(p-1)h,k+(q-1)h}$ for each $p,q=1,\ldots,d$ and $i,k=1,\ldots,h$.\newline\newline
We can now split the considered sum in \eqref{eq: break down original process} in a way such that we are able to express it in terms of sample cross-covariances afterwards. In order to do so, we define the sample cross-covariance at lag $l$ by 
\begin{align*}
\hat{r}^{(p,q)}_n(l):=\frac{1}{n}\sum_{j=1}^{n-l} X_j^{(p)}X_{j+l}^{(q)}
\end{align*}
 for $p,q=1,\ldots,d$.\newline
Note that in the case $h=1$, it follows directly that 
\begin{align*}
\sum_{j=1}^n \sum_{p,q=1}^d \sum_{i,k=1}^h Y_{j+i-1}^{(p)} Y_{j+k-1}^{(q)} \alpha_{ik}^{(p,q)}=\sum_{p,q=1}^d \alpha_{1,1}^{(p,q)}\sum_{j=1}^n  Y_j^{(p)}Y_j^{(q)}=n\sum_{p,q=1}^d  \hat{r}^{(p,q)}_n(0).
\end{align*}
The case $h=2$ has to be regarded separately, too, and we obtain
\begin{align*}
&\sum_{j=1}^n \sum_{p,q=1}^d \sum_{i,k=1}^2 Y_{j+i-1}^{(p)} Y_{j+k-1}^{(q)} \alpha_{ik}^{(p,q)}\\
=&\sum_{p,q=1}^d\left( \alpha_{1,1}^{(p,q)}\sum_{j=1}^n  Y_j^{(p)}Y_j^{(q)}+\alpha_{1,2}^{(p,q)}\sum_{j=1}^n  Y_j^{(p)}Y_{j+1}^{(q)}+\alpha_{2,1}^{(p,q)}\sum_{j=1}^n  Y_{j+1}^{(p)}Y_{j}^{(q)}+\alpha_{2,2}^{(p,q)}\sum_{j=1}^n Y_{j+1}^{(p)}Y_{j+1}^{(q)}\right)\\
=&\sum_{p,q=1}^d\Bigg( \alpha_{1,1}^{(p,q)} n \hat{r}_n^{(p,q)}(0)+\alpha_{1,2}^{(p,q)}\left(n\hat{r}_n^{(p,q)}(1)+\underbrace{Y_n^{(p)}Y_{n+1}^{(q)}}_{\bigstar}\right)+\alpha_{2,1}^{(p,q)}\left(n\hat{r}_n^{(q,p)}(1)+\underbrace{Y_{n+1}^{(p)}Y_{n}^{(q)}}_{\bigstar}\right)\\
&\quad\quad+\alpha_{2,2}^{(p,q)}\left(n \hat{r}_n^{(p,q)}(0)+\underbrace{Y_{n+1}^{(p)}Y_{n+1}^{(q)}}_{\bigstar}-\underbrace{Y_{1}^{(p)}Y_{1}^{(q)}}_{\bigstar}\right)\Bigg),
\end{align*}
Note that for each of the terms labeled by $\bigstar$ the following holds for $d^*\in\left(\frac{1}{4},\frac{1}{2}\right)$:
\begin{align*}
n^{-2d^*}\bigstar \xrightarrow{\mathbb{P}} 0,~(n\rightarrow\infty).
\end{align*}
We use this property later on when dealing with the asymptotics of the term in \eqref{eq: break down original process}.\newline
Finally we consider the term in \eqref{eq: break down original process} for $h\geq 3$ and arrive at
\begin{align}
&\sum_{j=1}^n \sum_{p,q=1}^d \sum_{i,k=1}^h Y_{j+i-1}^{(p)} Y_{j+k-1}^{(q)} \alpha_{ik}^{(p,q)}\notag\\
=&\sum_{p,q=1}^d \sum_{i,k=1}^h \alpha_{ik}^{(p,q)} \sum_{j=i}^{n+i-1} Y_{j}^{(p)} Y_{j+k-i}^{(q)} \notag\\
=&\sum_{p,q=1}^d \sum_{l=0}^{h-1} \sum_{i=1}^{h-l} \alpha_{i,i+l}^{(p,q)}\sum_{j=i}^{n+i-1} Y_j^{(p)}Y_{j+l}^{(q)}\notag\\
&\quad\quad+\sum_{p,q=1}^d   \sum_{l=-(h-1)}^{-1} \sum_{i=1-l}^{h} \alpha_{i,i+l}^{(p,q)}\sum_{j=i}^{n+i-1} Y_j^{(p)}Y_{j+l}^{(q)} \notag\\
=&\sum_{p,q=1}^d \sum_{l=0}^{h-1} \sum_{i=1}^{h-l} \alpha_{i,i+l}^{(p,q)}\sum_{j=i}^{n+i-1} Y_j^{(p)}Y_{j+l}^{(q)}\notag\\
&\quad\quad+\sum_{p,q=1}^d   \sum_{l=1}^{h-1} \sum_{i=1}^{h-l} \alpha_{i+l,i}^{(p,q)}\sum_{j=i}^{n+i-1} Y_{j+l}^{(p)}Y_{j}^{(q)} \label{eq: nochschön} \\
=&\sum_{p,q=1}^d \sum_{i=1}^h \alpha_{i,i}^{(p,q)} \sum_{j=i}^{n+i-1} Y_j^{(p)}Y_j^{(q)}\notag\\
&\quad\quad+\sum_{p,q=1}^d   \sum_{l=1}^{h-1} \sum_{i=1}^{h-l} \Bigg(\alpha_{i,i+l}^{(p,q)}\sum_{j=i}^{n+i-1} Y_j^{(p)}Y_{j+l}^{(q)}+ \alpha_{i+l,i}^{(p,q)}\sum_{j=i}^{n+i-1} Y_{j+l}^{(p)}Y_{j}^{(q)}  \Bigg)\label{eq: technical-l}\\
=&\sum_{p,q=1}^d \Bigg( \alpha_{1,1}^{(p,q)}\sum_{j=1}^n Y_j^{(p)}Y_j^{(q)}+\sum_{i=2}^h \alpha_{i,i}^{(p,q)} \sum_{j=i}^{n+i-1} Y_j^{(p)}Y_j^{(q)} \Bigg)\notag\\
&\quad\quad+\sum_{p,q=1}^d   \sum_{l=1}^{h-2} \Bigg(\Bigg( \alpha_{1,1+l}^{(p,q)}\sum_{j=1}^{n} Y_j^{(p)}Y_{j+l}^{(q)}+ \alpha_{1+l,1}^{(p,q)}\sum_{j=1}^{n} Y_{j+l}^{(p)}Y_{j}^{(q)}\Bigg) \notag\\
&\quad\quad\quad\quad\quad\quad+\sum_{i=2}^{h-l} \Bigg(\alpha_{i,i+l}^{(p,q)}\sum_{j=i}^{n+i-1} Y_j^{(p)}Y_{j+l}^{(q)}+ \alpha_{i+l,i}^{(p,q)}\sum_{j=i}^{n+i-1} Y_{j+l}^{(p)}Y_{j}^{(q)}  \Bigg)\Bigg)\notag\\
&\quad\quad+\sum_{p,q=1}^d \Bigg(\alpha_{1,h}^{(p,q)}\sum_{j=1}^n Y_j^{(p)}Y_{j+h-1}^{(q)}+\alpha_{h,1}^{(p,q)}Y_{j+h-1}^{(p)}Y_{j}^{(q)}\Bigg)
  \label{eq: technical-i}\\
=&\sum_{p,q=1}^d \Bigg( \alpha_{1,1}^{(p,q)} n \hat{r}_n^{(p,q)}(0)+\sum_{i=2}^h \alpha_{i,i}^{(p,q)} \Bigg(\underbrace{\sum_{j=n+1}^{n+i-1} Y_j^{(p)}Y_j^{(q)}}_{\bigstar} + n \hat{r}_n^{(p,q)}(0)-\underbrace{\sum_{j=1}^{i-1} Y_j^{(p)}Y_j^{(q)} }_{\bigstar} \Bigg)\Bigg)\notag\\
&+\sum_{p,q=1}^d   \sum_{l=1}^{h-2} \Bigg(\Bigg( \alpha_{1,1+l}^{(p,q)}n\hat{r}_n^{(p,q)}(l)+ \alpha_{1+l,1}^{(p,q)}n\hat{r}_n^{(q,p)}(l)\Bigg) \notag\\
&\quad\quad\quad\quad\quad\quad+\sum_{i=2}^{h-l} \Bigg(\alpha_{i,i+l}^{(p,q)}\Bigg(\underbrace{\sum_{j=n-l+1}^{n+i-1} Y_j^{(p)}Y_{j+l}^{(q)}}_{\bigstar}+n\hat{r}_n^{(p,q)}(l)-\underbrace{\sum_{j=1}^{i-1}Y_j^{(p)}Y_{j+l}^{(q)}}_{\bigstar}\Bigg)\notag\\
&\quad\quad\quad\quad\quad\quad+ \alpha_{i+l,i}^{(p,q)}\Bigg(\underbrace{\sum_{j=n-l+1}^{n+i-1} Y_{j+l}^{(p)}Y_{j}^{(q)}}_{\bigstar}+n\hat{r}_n^{(q,p)}(l)-\underbrace{\sum_{j=1}^{i-1}Y_{j+l}^{(p)}Y_{j}^{(q)}}_{\bigstar}\Bigg)  \Bigg)\Bigg)\notag\\
&+\sum_{p,q=1}^d \Bigg( \alpha_{1,h}^{(p,q)} \Bigg(\underbrace{\sum_{j=n-h+2}^{n} Y_{j}^{(p)}Y_{j+h-1}^{(q)}}_{\bigstar}+n \hat{r}_n^{(p,q)}(h-1)\Bigg)\\
&\quad\quad\quad\quad\quad\quad+\alpha_{h,1}^{(p,q)}\Bigg(\underbrace{\sum_{j=n-h+2}^{n} Y_{j+h-1}^{(p)}Y_{j}^{(q)}}_{\bigstar}+n\hat{r}_n^{(q,p)}(h-1)\Bigg) \Bigg).\label{eq: hässlichsteumformungever}
\end{align}
\noindent Again for each of the terms labeled by $\bigstar$ it holds for $d^*\in\left(\frac{1}{4},\frac{1}{2}\right)$:
\begin{align*}
n^{-2d^*}\bigstar \xrightarrow{\mathbb{P}} 0,~(n\rightarrow\infty),
\end{align*}
since each $\bigstar$ describes a sum with a finite number (independent of $n$) of summands. Therefore, we continue to express the terms denoted by $\bigstar$ by $o_{\mathbb{P}}\left(n^{2d^*}\right)$.\newline
With these calculations we are able to re-express the partial sum, whose asymptotics we are interested in, in terms of the sample cross-correlations of the original long-range dependent process $\left(Y_j\right)_{j\in\mathbb{Z}}$.  \newline
Finally, the previous calculations lead to
\begin{align}
&\sum_{j=1}^n f\left(Y_{j,h}\right)\notag\\
\overset{\mathcal{D}}{=}\:\:\,&\sum_{j=1}^n f^*\left(U_{j,h}\right)\notag\\
\overset{\eqref{eq: th2 umformung1}}{=}&\sum_{j=1}^n\left(\sum_{1\leq i,k\leq dh} \mathbb{E}\left( f^*\left(U_{j,h}\right) U_{j,h}^{(i)}  U_{j,h}^{(k)} \right)  U_{j,h}^{(i)}  U_{j,h}^{(k)} -\sum_{i=1}^{dh} \mathbb{E}\left(f^*\left(U_{j,h}\right) \left(U_{j,h}^{(i)}\right)^2\right)\right)+o_{\mathbb{P}}(n^{2d^*})\notag\\
\overset{\mathcal{D}}{\underset{\eqref{eq: break down original process}}{=}}&\sum_{j=1}^n \sum_{p,q=1}^d \sum_{i,k=1}^h  \alpha_{ik}^{(p,q)} \left(Y_{j+i-1}^{(p)} Y_{j+k-1}^{(q)}-\mathbb{E}\left(Y_{j+i-1}^{(p)} Y_{j+k-1}^{(q)}\right)\right)+o_{\mathbb{P}}(n^{2d^*}),\label{th1: umformung2}
\end{align}
where \eqref{th1: umformung2} follows, since \eqref{th1:erwartungswertgleichheit} yields
\begin{align*}
\sum_{i=1}^{dh} \mathbb{E}\left(f^*\left(U_{j,h}\right) \left(U_{j,h}^{(i)}\right)^2\right)=\:\:\,&\sum_{1\leq i,k\leq dh} \mathbb{E}\left( f^*\left(U_{j,h}\right) U_{j,h}^{(i)}  U_{j,h}^{(k)} \right)  \mathbb{E}\left(U_{j,h}^{(i)}  U_{j,h}^{(k)}\right)\\
\overset{\eqref{eq: break down original process}}{=}&\sum_{p,q=1}^d \sum_{i,k=1}^h  \alpha_{ik}^{(p,q)}\mathbb{E}\left(Y_{j+i-1}^{(p)} Y_{j+k-1}^{(q)}\right).
\end{align*}
Taking the parts containing the sample cross-correlations into account, we derive
\begin{align}
&\sum_{j=1}^n \sum_{p,q=1}^d \sum_{i,k=1}^h  \alpha_{ik}^{(p,q)} \left(Y_{j+i-1}^{(p)} Y_{j+k-1}^{(q)}-\mathbb{E}\left(Y_{j+i-1}^{(p)} Y_{j+k-1}^{(q)}\right)\right)+o_{\mathbb{P}}(n^{2d^*})\notag\\
\overset{\eqref{eq: hässlichsteumformungever}}{=}&\sum_{p,q=1}^d \left( \alpha_{1,1}^{(p,q)} n \left(\hat{r}_n^{(p,q)}(0)-r^{(p,q)}(0)\right)+\sum_{i=2}^h \alpha_{i,i}^{(p,q)}  n  \left(\hat{r}_n^{(p,q)}(0)-r^{(p,q)}(0)\right)\right)\notag\\
&+\sum_{p,q=1}^d   \sum_{l=1}^{h-2} \Bigg(\left( \alpha_{1,1+l}^{(p,q)}n \left(\hat{r}_n^{(p,q)}(l)-r^{(p,q)}(l)\right)+ \alpha_{1+l,1}^{(p,q)}n \left(\hat{r}_n^{(q,p)}(l)-r^{(q,p)}(l)\right)\right) \notag\\
&\quad\quad\quad\quad\quad\quad+\sum_{i=2}^{h-l} \Bigg(\alpha_{i,i+l}^{(p,q)}n\left(\hat{r}_n^{(p,q)}(l)-r^{(p,q)}(l)\right)+ \alpha_{i+l,i}^{(p,q)} n\left(\hat{r}_n^{(q,p)}(l)-r^{(q,p)}(l)\right)  \Bigg)\Bigg)\notag\\
&+\sum_{p,q=1}^d \left( \alpha_{1,h}^{(p,q)}n\left(\hat{r}_n^{(p,q)}(h-1)-r^{(p,q)}(h-1)\right)+\alpha_{h,1}^{(p,q)} n\left(\hat{r}_n^{(q,p)}(h-1)-r^{(q,p)}(h-1)\right) \right)\\
&+o_{\mathbb{P}}(n^{2d^*})\notag\\
\notag\\
=\:\:\,&n\sum_{p,q=1}^d \Bigg(\sum_{l=0}^{h-1}\sum_{i=1}^{h-l}  \alpha_{i,i+l}^{(p,q)} \left(\hat{r}_n^{(p,q)}(l)-r^{(p,q)}(l)\right)+  \sum_{l=1}^{h-1}\sum_{i=1}^{h-l} \alpha_{i+l,i}^{(p,q)} \left(\hat{r}_n^{(q,p)}(l)-r^{(q,p)}(l)\right)  \Bigg)\Bigg)\notag\\
&\quad+o_{\mathbb{P}}(n^{2d^*}).\label{eqth2:1}
\end{align}
We take a closer look at the impact of each long-range dependence parameter $d_p$, $p=1,\ldots,d$ to the convergence of this sum. The setting we are considering does not allow for a normalization depending on $p$ and $q$ for each cross-correlation $\left(\hat{r}_n^{(p,q)}(l)-r^{(p,q)}(l)\right)$, $l=0,\ldots,h-1$ but we need to find a normalization for all $p,q=1,\ldots, d$. Hence, we need to remember the set $P^*:=\{p\in\{1,\ldots,d\}:~d_p\geq d_q \forall q\in\{1,\ldots,d\}\}$ and the parameter $d^*=\max\limits_{p=1,\ldots,d} d_p$, such that for each $p\in P^*$ we have $d_p=d^*$. For each $p,q\in\{1,\ldots,d\}$ with $(p,q)\notin P^*\times P^*$ and $l=0,\ldots,h-1$, we conclude that
\begin{align}
\mathbb{E}\left(\left( n^{1-2d^*}\left(\hat{r}_n^{(p,q)}(l)-r^{(p,q)}(l)\right)\right)^2\right)&=n^{2\left(d_p+d_q-2d^*\right)}\mathbb{E}\left( \left( n^{1-d_p-d_q}\left(\hat{r}_n^{(p,q)}(l)-r^{(p,q)}(l)\right)\right)^2\right)\notag\\
&= n^{2d_p+2d_q-4d^*} C_2 \left(L_{p,p}L_{q,q}+L_{p,q}L_{q,p}\right)\notag\\
&\xrightarrow {(n\rightarrow\infty)}0\label{convergence0},  
\end{align}
since $d_p+d_q-2d^*<0$.\newline
This implies that 
\begin{align*}
 n^{1-2d^*}\left(\hat{r}_n^{(p,q)}(0)-r^{(p,q)}(0)\right)\xrightarrow{\mathbb{P}} 0
\end{align*}
Hence, using Slutsky's theorem, the crucial parameters that determine the normalization and, therefore, the limit distribution of \eqref{eq:bisalphaschlange} are given in $P^*$. We have an equal long-range dependence parameter $d^*$ to regard for all $p\in P^*$. Applying Lemma \ref{lemm:HR2preceding}, we obtain the following, by using the symmetry in $l=0$ of the cross correlation function $r^{(p,q)}(0)=r^{(q,p)}(0)$ for $p,q\in P^*$:
\begin{align}
&\sum_{p,q=1}^d \Bigg(\sum_{l=0}^{h-1}\sum_{i=1}^{h-l}  \alpha_{i,i+l}^{(p,q)}\left(\hat{r}_n^{(p,q)}(l)-r^{(p,q)}(l)\right)+  \sum_{l=1}^{h-1}\sum_{i=1}^{h-l} \alpha_{i+l,i}^{(p,q)} \left(\hat{r}_n^{(q,p)}(l)-r^{(q,p)}(l)\right)  \Bigg)\Bigg)\notag\\
=& \sum_{p,q\in P^*} \Bigg(\sum_{l=0}^{h-1}\sum_{i=1}^{h-l}  \alpha_{i,i+l}^{(p,q)}\left(\hat{r}_n^{(p,q)}(0)-r^{(p,q)}(0)\right)+  \sum_{l=1}^{h-1}\sum_{i=1}^{h-l} \alpha_{i+l,i}^{(p,q)} \left(\hat{r}_n^{(q,p)}(0)-r^{(q,p)}(0)\right)  \Bigg)\Bigg)+o_{\mathbb{P}}(n^{2d^*-1})\notag\\
=& \sum_{p,q\in P^*}\left(\hat{r}_n^{(p,q)}(0)-r^{(p,q)}(0)\right) \Bigg(\sum_{l=0}^{h-1}\sum_{i=1}^{h-l}  \alpha_{i,i+l}^{(p,q)}+  \sum_{l=1}^{h-1}\sum_{i=1}^{h-l} \alpha_{i+l,i}^{(p,q)} \Bigg)+o_{\mathbb{P}}(n^{2d^*})\notag\\
=& \sum_{p,q\in P^*}\left(\hat{r}_n^{(p,q)}(0)-r^{(p,q)}(0)\right) \left(\sum_{i,k=1}^h \alpha_{i,k}^{(p,q)}\right)+o_{\mathbb{P}}(n^{2d^*})\notag\\
=& \sum_{p,q\in P^*}\tilde{\alpha}^{(p,q)}\left(\hat{r}_n^{(p,q)}(0)-r^{(p,q)}(0)\right)+o_{\mathbb{P}}(n^{2d^*}), \label{eq:bisalphaschlange}
\end{align}
 by defining $\tilde{\alpha}^{(p,q)}:=\sum\limits_{i,k=1}^h \alpha_{i,k}^{(p,q)}$. 
\noindent Applying the continuous mapping theorem given in \cite{vandervaart:2000}, Theorem 2.3 to the result in Corollary \ref{korl:precedingHR2} we arrive at
\begin{align*}
n^{-2d^*}\left(C_2\right)^{-1/2}\sum_{j=1}^n f\left(Y_{j,h}\right)&= n^{-2d^*} \left(n \sum_{p,q=1}^d\tilde{\alpha}^{(p,q)}\left(\hat{r}_n^{(p,q)}(0)-r^{(p,q)}(0)\right) +o_{\mathbb{P}}(n^{2d^*})\right)\\
&=n^{1-2d^*}\left(C_2\right)^{-1/2} \sum_{p,q=1}^d\tilde{\alpha}^{(p,q)}\left(\hat{r}_n^{(p,q)}(0)-r^{(p,q)}(0)\right) + o_{\mathbb{P}}(1)\notag\\
&\xrightarrow{\mathcal{D}}\sum_{p,q\in P^*}\tilde{\alpha}^{(p,q)} Z^{(p,q)}_{2,d^*+1/2}(1),
\end{align*}
where 
\begin{align*}
Z^{(p,q)}_{2,d^*+1/2}(1)=K_{p,q}\left(d^*\right)\int_{\mathbb{R}^2}^{\prime\prime} \frac{\exp\left(i\left(\lambda_1+\lambda_2\right)\right)-1}{i\left(\lambda_1+\lambda_2\right)}\left|\lambda_1\lambda_2  \right|^{-d^*} \tilde{B}_L^{(p)}\left(\mathrm{d}\lambda_1\right)\tilde{B}_L^{(q)}\left( \mathrm{d}\lambda_2\right).
\end{align*}
The matrix $K\left(d^*\right)$ is given in Corollary \ref{korl:precedingHR2}. Moreover, $\tilde{B}_L(\mathrm{d}\lambda)$ is a multivariate Hermitian-Gaussian random measure with $\mathbb{E}\left(B_L(\mathrm{d}\lambda)B_L(\mathrm{d}\lambda)^*\right)=L~\mathrm{d}\lambda$ and $L$ as defined in (2).
$\hfill\Box$\newline
\end{pf}
    \renewcommand\thesection{C}%
\section{Proof of Corollary 2.2}
\begin{pf}
We assumed $d^*\in\left(\frac{1}{4},\frac{1}{2}\right)$, because otherwise we leave the long-range dependent setting, since we are studying functionals with Hermite rank $2$ and the transformed process would no longer be long-range dependent and limit theorems for functionals of short-range dependent processes would hold, see Theorem 4 in \cite{arcones:1994}. This choice of $d^*$ assures that the multivariate generalization of the Reduction theorem as it is used in the proof of Theorem 2.1 still holds for these softened assumptions, as explained in (8).\newline
We turn to the asymptotics of $g^{(p,q)}\left(Y_j\right)$. We obtain for all $p,q\in\{1,\ldots,d\}\setminus P^*$,i.e., excluding $d_p=d_q=d^*$ and for all $l=0,\ldots,h-1$ as in \eqref{convergence0}, that
\begin{align}
\mathbb{E}\left(\left( n^{1-2d^*}\left(\hat{r}_n^{(p,q)}(l)-r^{(p,q)}(l)\right)\right)^2\right)&=n^{2\left(d_p+d_q-2d^*\right)}\mathbb{E}\left( \left( n^{1-d_p-d_q}\left(\hat{r}_n^{(p,q)}(l)-r^{(p,q)}(l)\right)\right)^2\right)\notag\\
&= n^{2d_p+2d_q-4d^*} C_2 \left(L_{p,p}L_{q,q}+L_{p,q}L_{q,p}\right)\notag\\
&\xrightarrow {(n\rightarrow\infty)}0,\label{eq:convergence1}
\end{align}
since $d_p+d_q-2d^*<0$.\newline
This implies that 
\begin{align*}
 n^{1-2d^*}\left(\hat{r}_n^{(p,q)}(0)-r^{(p,q)}(0)\right)\xrightarrow{\mathbb{P}} 0.
\end{align*}
Applying Slutsky's theorem, we observe that only $p,q \in P^*$ have an impact on the convergence behaviour as it is given in \eqref{eq:bisalphaschlange} and hence, the result in Theorem 2.1 holds.
 $\hfill\Box$\newline
\end{pf}
    \renewcommand\thesection{D}%
\section{Proof of Theorem 2.3}
\begin{pf}
We follow the proof of Theorem 2.1 until \eqref{eq:bisalphaschlange}, in order to obtain a limit distribution that can be expressed by the sum of two standard Rosenblatt random variables: 
\begin{align}
&\sum_{p,q=1}^2 \tilde{\alpha}^{(p,q)} \left(\hat{r}_n^{(p,q)}(0)-r^{(p,q)}(0)\right)\notag\\
=&\frac{1}{n}\sum_{j=1}^n\sum_{p,q=1}^2 \tilde{\alpha}^{(p,q)} \left(Y_j^{(p)}Y_j^{(q)}-r^{(p,q)}(0)\right)\notag\\
=&\frac{1}{n}\sum_{j=1}^n \left(Y_j^{(1)},Y_j^{(2)}\right) \begin{pmatrix}
 \tilde{\alpha}^{(1,1)} && \tilde{\alpha}^{(1,2)}  \\ \tilde{\alpha}^{(2,1)}  && \tilde{\alpha}^{(2,2)} 
\end{pmatrix}  \left(Y_j^{(1)},Y_j^{(2)}\right)^t\notag\\
&\quad\quad\quad\quad\quad\quad\quad\quad\quad-\mathbb{E}\left(\left(Y_j^{(1)},Y_j^{(2)}\right) \begin{pmatrix}
 \tilde{\alpha}^{(1,1)} && \tilde{\alpha}^{(1,2)}  \\ \tilde{\alpha}^{(2,1)}  && \tilde{\alpha}^{(2,2)} 
\end{pmatrix} \left(Y_j^{(1)},Y_j^{(2)}\right)^t\right). \label{eq: matrixumformung}
\end{align}
We remember that $\tilde{\alpha}^{(p,q)}=\sum_{i,k=1}^h \alpha_{i,k}^{(p,q)}=\sum_{i,k=1}^h \alpha_{i+(p-1)h,k+(q-1)h}$ for $p,q=1,2$ and $\mathbb{A}=\left(\alpha_{i,k}\right)_{1\leq i,k\leq 2h}=\Sigma_{2,h}^{-1}C\Sigma_{2,h}^{-1}$. Since $\Sigma_{2,h}^{-1}$ is the inverse of the covariance matrix $\Sigma_{2,h}$ of $Y_{1,h}$ it is a symmetric matrix. The matrix of second order Hermite coefficients $C$ has the representation $C=\mathbb{E}\left(Y_{j,h} f\left(Y_{j,h} \right)Y_{j,h}^t\right)$ and, therefore, $c_{i,k}=\mathbb{E}\left(Y_{j,h}^{(i)}Y_{j,h}^{(k)}f\left(Y_{j,h}\right)\right)=c_{k,i}$ for each $i,k=1,\ldots,2h$.
Then, $\mathbb{A}$ is a symmetric matrix, too, since $\mathbb{A}^t=\left(\Sigma_{2,h}^{-1}C\Sigma_{2,h}^{-1}\right)^t=\left(\Sigma_{2,h}^{-1}\right)^tC^t\left(\Sigma_{2,h}^{-1}\right)^t=\mathbb{A}$. We can now show that $ \begin{pmatrix}
 \tilde{\alpha}^{(1,1)} && \tilde{\alpha}^{(1,2)}  \\ \tilde{\alpha}^{(2,1)}  && \tilde{\alpha}^{(2,2)} 
\end{pmatrix}$ is a symmetric matrix, i.e., $ \tilde{\alpha}^{(1,2)}= \tilde{\alpha}^{(2,1)}$. To this end, we define $\mathbb{I}_p=\left(0,0,\ldots,0,1,\ldots,1,0,\ldots,0\right)^t\in\mathbb{R}^{2h}$ such that $\mathbb{I}_p^{(i)}=1$ only if $i=(p-1)h+1,\ldots,ph$, $p=1,2$. Then, we obtain
\begin{align*}
\tilde{\alpha}^{(1,2)}=\sum_{i,k=1}^h \alpha_{i,k}^{(1,2)}=\left(\tilde{\alpha}^{(1,2)}\right)^t=\left(\mathbb{I}_1^t \mathbb{A} \mathbb{I}_2\right)^t=\mathbb{I}_2^t \mathbb{A}\mathbb{I}_1=\tilde{\alpha}^{(2,1)}.
\end{align*}
We now apply the new assumption that $r^{(1,1)}(l)=r^{(2,2)}(l)$, for $l=0,\ldots,h-1$ and show $\tilde{\alpha}^{(1,1)}=\tilde{\alpha}^{(2,2)}$ with the symmetry features of the multivariate normal distribution discussed in (2.2) and in (2.3) in \cite{nuessgen:2021}, since
$c_{i,j}=c_{2h-i+1,2h-j+1}$, $i,j=1,\ldots,2h$.\newline
We have to study
\begin{align*}
\tilde{\alpha}^{(2,2)}=\left(\mathbb{I}_2^t \mathbb{A} \mathbb{I}_2\right)^t=\mathbb{I}_2^t \Sigma_{2,h}^{-1}C\Sigma_{2,h}^{-1} \mathbb{I}_2.
\end{align*}
Since $\Sigma_{2,h}^{-1}=\left(g_{i,k}\right)_{1\leq i,k\leq 2h}$ is a symmetric and persymmetric matrix, we have $g_{i,k}=g_{k,i}$ and $g_{i,k}=g_{2h-i+1,2h-k+1}$ for $i,k=1,\ldots,2h$. Then, we obtain
\begin{align*}
\mathbb{I}_2^t\Sigma_{2,h}^{-1}&=\left(\sum_{i=h+1}^{2h} g_{i,1},\ldots,\sum_{i=h+1}^{2h} g_{i,2h}\right)\\
&=\left(\sum_{i=1}^{h} g_{i+h,1},\ldots,\sum_{i=1}^{h} g_{i+h,2h}\right)\\
&=\left(\sum_{i=1}^{h} g_{h-i+1,2h},\ldots,\sum_{i=1}^{h} g_{h-i+1,1}\right)\\
&=\left(\sum_{i=1}^{h} g_{i,2h},\ldots,\sum_{i=1}^{h} g_{i,1}\right)\\
&=\left(\sum_{i=1}^{h} g_{2h,i},\ldots,\sum_{i=1}^{h} g_{1,i}\right)\\
&=:\left(\tilde{g}_{2h},\ldots,\tilde{g}_1\right).
\end{align*}
Note that
\begin{align*}
\Sigma_{2,h}^{-1}\mathbb{I}_1=\left(\sum_{i=1}^h g_{1,i},\ldots,\sum_{i=1}^h g_{2h,i}\right)^t=\left(\tilde{g}_{1},\ldots,\tilde{g}_{2h}\right)^t.
\end{align*}
Then, we arrive at
\begin{align*}
\tilde{\alpha}^{(2,2)}=\left(\mathbb{I}_2^t \mathbb{A} \mathbb{I}_2\right)^t&=\mathbb{I}_2^t \Sigma_{2,h}^{-1}C\Sigma_{2,h}^{-1} \mathbb{I}_2\\
&=\sum_{i,k=1}^{2h} \tilde{g}_{2h-i+1}\tilde{g}_{2h-k+1} c_{i,k}\\
&=\sum_{i,k=1}^{2h} \tilde{g}_{2h-i+1}\tilde{g}_{2h-k+1} c_{2h-i+1,2h-k+1}\\
&=\sum_{i,k=1}^{2h} \tilde{g}_{i}\tilde{g}_{k} c_{i,k}\\
&=\mathbb{I}_1^t \Sigma_{2,h}^{-1}C\Sigma_{2,h}^{-1} \mathbb{I}_1\\
&=\tilde{\alpha}^{(1,1)}.
\end{align*}
So we have to deal with a special type of $2\times2$-matrix, since the original matrix in the formula \eqref{eq:bisalphaschlange}, namely $ \begin{pmatrix}
 \tilde{\alpha}^{(1,1)} && \tilde{\alpha}^{(1,2)}  \\ \tilde{\alpha}^{(2,1)}  && \tilde{\alpha}^{(2,2)}
\end{pmatrix}$ has now reduced to $ \begin{pmatrix}
 \tilde{\alpha}^{(1,1)} && \tilde{\alpha}^{(1,2)}  \\ \tilde{\alpha}^{(1,2)}  && \tilde{\alpha}^{(1,1)} 
\end{pmatrix}$.\newline
Finally, we know that any real-valued symmetric matrix $A$ can be decomposed via diagonalization into an orthogonal matrix $V$ and a diagonal matrix $D$, where the entries of the latter one are determined via the eigenvalues of $A$, for details, see \cite{beutelspacher:2003}, p. 327.\newline
We can explicity give formulas for the entries of these matrices here:
\begin{align*}
V= \begin{pmatrix}
-2^{-1/2} && 2^{-1/2}  \\ 2^{-1/2}  && 2^{-1/2}
\end{pmatrix},\quad
D=\begin{pmatrix}
\lambda_1=\tilde{\alpha}^{(1,1)}-\tilde{\alpha}^{(1,2)} && 0 \\ 0 && \lambda_2=\tilde{\alpha}^{(1,1)}+\tilde{\alpha}^{(1,2)}
\end{pmatrix},
\end{align*}
such that
\begin{align*}
VDV=\begin{pmatrix}
\frac{\lambda_1+\lambda_2}{2} && \frac{\lambda_2-\lambda_1}{2} \\
\frac{\lambda_2-\lambda_1}{2} && \frac{\lambda_1+\lambda_2}{2}
\end{pmatrix}=\begin{pmatrix}
 \tilde{\alpha}^{(1,1)} && \tilde{\alpha}^{(1,2)}  \\ \tilde{\alpha}^{(1,2)}  && \tilde{\alpha}^{(1,1)} 
\end{pmatrix}.
\end{align*}
So continuing with \eqref{eq: matrixumformung}, we now have the representation
\begin{align}
&\frac{1}{n}\sum_{j=1}^n \left(Y_j^{(1)},Y_j^{(2)}\right) \begin{pmatrix}
 \tilde{\alpha}^{(1,1)} && \tilde{\alpha}^{(1,2)}  \\ \tilde{\alpha}^{(2,1)}  && \tilde{\alpha}^{(2,2)} 
\end{pmatrix}  \left(Y_j^{(1)},Y_j^{(2)}\right)^t\notag\\
&\quad\quad\quad\quad\quad\quad\quad\quad\quad-\mathbb{E}\left(\left(Y_j^{(1)},Y_j^{(2)}\right) \begin{pmatrix}
 \tilde{\alpha}^{(1,1)} && \tilde{\alpha}^{(1,2)}  \\ \tilde{\alpha}^{(2,1)}  && \tilde{\alpha}^{(2,2)} 
\end{pmatrix} \left(Y_j^{(1)},Y_j^{(2)}\right)^t\right)\notag\\
=&\frac{1}{n}\sum_{j=1}^n \left(Y_j^{(1)},Y_j^{(2)}\right) VDV  \left(Y_j^{(1)},Y_j^{(2)}\right)^t-\mathbb{E}\left(\left(Y_j^{(1)},Y_j^{(2)}\right)VDV\left(Y_j^{(1)},Y_j^{(2)}\right)^t\right)\notag\\
=&\frac{1}{n}\sum_{j=1}^n \frac{\tilde{\alpha}^{(1,1)}-\tilde{\alpha}^{(1,2)}}{2} \left(\left(Y_j^{(2)}-Y_j^{(1)}\right)^2-\mathbb{E}\left(Y_j^{(2)}-Y_j^{(1)}\right)^2\right)\notag\\
&\quad\quad +\frac{1}{n}\sum_{j=1}^n\frac{\tilde{\alpha}^{(1,1)}+\tilde{\alpha}^{(1,2)}}{2}\left(\left(Y_j^{(1)}+Y_j^{(2)}\right)^2-\mathbb{E}\left(Y_j^{(1)}+Y_j^{(2)}\right)^2\right)\notag\\
=&\frac{1}{n}\sum_{j=1}^n\left(\tilde{\alpha}^{(1,1)}-\tilde{\alpha}^{(1,2)}\right)\left(1-r^{(1,2)}(0)\right) \left(\left(\frac{Y_j^{(2)}-Y_j^{(1)}}{\sqrt{2-2r^{(1,2)}(0)}}\right)^2-1\right)\notag\\
&\quad\quad +\frac{1}{n}\sum_{j=1}^n \left(\tilde{\alpha}^{(1,1)}+\tilde{\alpha}^{(1,2)}\right)\left(1+r^{(1,2)}(0)\right) \left(\left(\frac{Y_j^{(1)}+Y_j^{(2)}}{\sqrt{2+2r^{(1,2)}(0)}}\right)^2-1\right)\notag\\
&=\frac{1}{n}\left(\tilde{\alpha}^{(1,1)}-\tilde{\alpha}^{(1,2)}\right)\left(1-r^{(1,2)}(0)\right) \sum_{j=1}^n H_2\left(Y^*_j\right) \notag\\
&\quad\quad +\frac{1}{n}\left(\tilde{\alpha}^{(1,1)}+\tilde{\alpha}^{(1,2)}\right)\left(1+r^{(1,2)}(0)\right) \sum_{j=1}^n H_2\left(Y^{**}_j\right),\label{eq:alternativeprocesses}
\end{align}
with $Y^*_j:=\frac{Y_j^{(2)}-Y_j^{(1)}}{\sqrt{2-2r^{(1,2)}(0)}}$ and $Y^{**}_j:=\frac{Y_j^{(1)}+Y_j^{(2)}}{\sqrt{2+2r^{(1,2)}(0)}}$.\newline
Now note that 
\begin{align*}
\mathbb{E}\left(Y^*_j Y^{**}_j\right)=\mathbb{E}\left(\frac{Y_j^{(2)}-Y_j^{(1)}}{\sqrt{2-2r^{(1,2)}(0)}}\frac{Y_j^{(1)}+Y_j^{(2)}}{\sqrt{2+2r^{(1,2)}(0)}}\right)=0.
\end{align*}
Therefore, we created a bivariate long-range dependent Gaussian process, since
\begin{align*}
\begin{pmatrix}
-1 && 1 \\ 1 && 1 
\end{pmatrix}
\left(Y_j^{(1)}, Y_j^{(2)}\right)^t=\left(Y_j^*,Y_j^{**}\right)^t\sim \mathcal{N}\left(0,I_2\right)
\end{align*}
with cross-covariance function
\begin{align}
r_*^{(1,2)}(k):=\mathbb{E}\left(Y_j^*Y_{j+k}^{**}\right)&=\mathbb{E}\left(\frac{Y_j^{(2)}-Y_j^{(1)}}{\sqrt{2-2r^{(1,2)}(0)}}\frac{Y_{j+k}^{(1)}+Y_{j+k}^{(2)}}{\sqrt{2+2r^{(1,2)}(0)}}\right)\notag\\
&=\frac{r^{(2,1)}(k)+r^{(2,2)}(k)-r^{(1,1)}(k)-r^{(1,2)}(k)}{2\sqrt{\left(1-r^{(1,2)}(0)\right)\left(1+r^{(1,2)}(0)\right)}}\notag\\
&\simeq \frac{L_{2,2}+L_{2,1}-L_{1,2}-L_{1,1}}{2\sqrt{\left(1-r^{(1,2)}(0)\right)\left(1+r^{(1,2)}(0)\right)}}k^{2d^*-1}.\label{eq: alternativeCC}
\end{align}
Note that the covariance functions have the following asymptotic behaviour:
\begin{align*}
r_*^{(1,1)}(k):=\mathbb{E}\left(Y_j^*Y_{j+k}^{*}\right)&=\mathbb{E}\left(\frac{Y_j^{(2)}-Y_j^{(1)}}{\sqrt{2-2r^{(1,2)}(0)}}\frac{Y_{j+k}^{(2)}-Y_{j+k}^{(1)}}{\sqrt{2-2r^{(1,2)}(0)}}\right)\\
&=\frac{r^{(2,2)}(k)-r^{(2,1)}(k)-r^{(1,2)}(k)+r^{(1,1)}(k)}{2-2r^{(1,2)}(0)}\\
&\simeq \underbrace{ \frac{L_{2,2}-L_{2,1}-L_{1,2}+L_{1,1}}{2-2r^{(1,2)}(0)}}_{=:L_{1,1}^*}k^{2d^*-1}
\end{align*}
and analogously
\begin{align*}
r_*^{(2,2)}(k):=\mathbb{E}\left(Y_j^{**}Y_{j+k}^{**}\right)\simeq \underbrace{\frac{L_{2,2}+L_{2,1}+L_{1,2}+L_{1,1}}{2+2r^{(1,2)}(0)}}_{=:L_{2,2}^{*}}k^{2d^*-1}.
\end{align*}
We can now apply the result of \cite{arcones:1994}, Theorem 6, since we created a bivariate Gaussian process with independent entries for fixed $j$. Note that for the function we apply here, namely \newline
$\tilde{f}\left(Y_j^*,Y_j^{**}\right)=H_2\left(Y_j^{*}\right)+H_2\left(Y_j^{**}\right)$ the weighting factors in \cite{arcones:1994}, Theorem 6 reduce to $e_{1,1}=e_{2,2}=1$ and $e_{1,2}=e_{2,1}=0$. These weighting factors fit into the result in \cite{arcones:1994}, (3.6) and (3.7), that even yields joint convergence of the vector of both univariate summands, $\left(H_2\left(Y_j^{*}\right), H_2\left(Y_j^{**}\right)\right)$, suitably normalized, to a vector of two (dependent) Rosenblatt random variables. Since the long-range dependence property in \cite{kechagias:pipiras:2015}, Definition 2.1 is more specific than in \cite{arcones:1994}, p. 2259, (3.1) (see considerations in (8)), we are able to scale the variances of each Rosenblatt random variable to $1$ and give the covariance between them, by using the normalization given in \cite{beran:feng:ghosh:kulik:2013}, Theorem 4.3. We obtain
\begin{align*}
&n^{-2d*}\left(2C_2\right)^{-1/2}\left(\tilde{\alpha}^{(1,1)}-\tilde{\alpha}^{(1,2)}\right)\left(1-r^{(1,2)}(0)\right) \sum_{j=1}^n H_2\left(Y^*_j\right) \\
&\quad\quad +n^{-2d*}\left(2C_2\right)^{-1/2}\left(\tilde{\alpha}^{(1,1)}+\tilde{\alpha}^{(1,2)}\right)\left(1+r^{(1,2)}(0)\right) \sum_{j=1}^n H_2\left(Y^{**}_j\right)\\
&\xrightarrow{\mathcal{D}} \left(\tilde{\alpha}^{(1,1)}-\tilde{\alpha}^{(1,2)}\right)\left(1-r^{(1,2)}(0)\right)L_{1,1}^*Z^*_{2,d^*+1/2}(1)+\left(\tilde{\alpha}^{(1,1)}+\tilde{\alpha}^{(1,2)}\right)\left(1+r^{(1,2)}(0)\right) L_{2,2}^*Z^{**}_{2,d^*+1/2}\\
&=\left(\tilde{\alpha}^{(1,1)}-\tilde{\alpha}^{(1,2)}\right) \frac{L_{2,2}-L_{2,1}-L_{1,2}+L_{1,1}}{2}Z^*_{2,d^*+1/2}(1)\\
&\quad\quad+\left(\tilde{\alpha}^{(1,1)}+\tilde{\alpha}^{(1,2)}\right) \frac{L_{2,2}+L_{2,1}+L_{1,2}+L_{1,1}}{2}Z^{**}_{2,d^*+1/2}(1)
\end{align*}
with $C_2:=\frac{1}{2d^*\left(4d^*-1\right)}$ being the same normalizing factor as in Theorem 2.1.\newline
We observe that $Z^{*}_{2,d^*+1/2}(1)$ and $Z^{**}_{2,d^*+1/2}(1)$ are both standard Rosenblatt random variables. Following Corollary \ref{korl:precedingHR2}, their covariance is given by
\begin{align*}
\mathrm{Cov}\left(Z^{*}_{2,d^*+1/2}(1),Z^{**}_{2,d^*+1/2}(1)\right)&=\frac{\left(L_{1,2}^*+L_{2,1}^*\right)^2}{L_{1,1}^*L_{2,2}^*}\\
&=\frac{2\left( \left(L_{2,2}-L_{1,1}\right)^2  \right)}{4\left(1-r^{(1,2)}(0)\right)\left(1+r^{(1,2)}(0)\right)} \left(L_{1,1}^*L_{2,2}^*\right)^{-1} \\
&=\frac{\left(L_{2,2}-L_{1,1}\right)^2 }{\left(L_{1,1}+L_{2,2}\right)^2-\left(L_{1,2}+L_{2,1}\right)^2}.
\end{align*}
Note that $\left(L_{1,1}+L_{2,2}\right)^2-\left(L_{1,2}+L_{2,1}\right)^2\neq 0$ is fulfilled since $L_{1,1}+L_{2,2}\neq L_{1,2}+L_{2,1}$. 
$\hfill\Box$\newline
\end{pf}
\bibliographystyle{elsarticle-harv} 
\bibliography{OPDincontextLRD.bib}